\documentclass[11pt,a4paper]{article}
\usepackage{amssymb}
\usepackage{eurosym}
\usepackage{amsfonts}
\usepackage{amsmath}
\usepackage{amsthm}
\usepackage{graphicx}
\usepackage{float}
\usepackage{hyperref}

\hypersetup{
	colorlinks=true,
	linkcolor=blue,
	anchorcolor=blue,
	citecolor=blue
}
\newtheorem{theorem}{Theorem}[section]

\newtheorem{corollary}[theorem]{Corollary}
\newtheorem{lemma}[theorem]{Lemma}

\theoremstyle{definition}
\newtheorem{definition}[theorem]{Definition}

\newtheorem{remark}[theorem]{Remark}

\renewenvironment{proof}[1][Proof]{\noindent\textbf{#1.} }{\ \rule{0.5em}{0.5em}}

\renewcommand{\theequation}{\thesection.\arabic{equation}}
\allowdisplaybreaks
\input{tcilatex}

\def\limfunc#1{\mathop{\mathrm{#1}}}
\def\func#1{\mathop{\mathrm{#1}}\nolimits}

\def\Xint#1{\mathchoice
{\XXint\displaystyle\textstyle{#1}}%
{\XXint\textstyle\scriptstyle{#1}}%
{\XXint\scriptstyle\scriptscriptstyle{#1}}%
{\XXint\scriptscriptstyle\scriptscriptstyle{#1}}%
\!\int}
\def\XXint#1#2#3{{\setbox0=\hbox{$#1{#2#3}{\int}$ }
\vcenter{\hbox{$#2#3$ }}\kern-.6\wd0}}

\def\oint{\Xint-}

\def\enddoc{\end{document}}

\def\Qlb#1{#1}

\def\Qcb#1{#1} 

\def\FRAME#1#2#3#4#5#6#7#8
{
 \begin{figure}[H]
 \begin{center}
 \includegraphics[height=#3]{#7}
 \caption{#5}
 \label{#6}
 \end{center}
 \end{figure}
}

\begin{document}
	\title{Existence results for Leibenson's equation on Riemannian manifolds}
	\author{Philipp S\"urig}
	\date{June 2026}
	\maketitle
	
	\begin{abstract}
		We consider on an arbitrary Riemannian manifold $M$ the \textit{Leibenson equation} $\partial _{t}u=\Delta _{p}u^{q}$, that is also known as a \textit{doubly nonlinear evolution equation}. We prove that if $p>1$ and $q>0$ then the Cauchy-problem \begin{equation*}
			\left\{ 
			\begin{array}{ll}
				\partial _{t}u=\Delta _{p}u^{q} &\text{in}~M\times (0, \infty), \\u(x, 0)=u_{0}(x)& \text{in}~M,
			\end{array}%
			\right.  \end{equation*}
		has a unique weak solution for any  $u_{0}\in L^{1}(M)\cap L^{\infty}(M)$.
	\end{abstract}
	
	\let\thefootnote\relax\footnotetext{
		The author was funded by the Deutsche Forschungsgemeinschaft (DFG,
		German Research Foundation) - Project-ID 317210226 - SFB 1283.}
	
	\tableofcontents
	
\section{Introduction}
Let $M$ be an arbitrary geodesically complete Riemannian manifold. We are concerned here with a non-linear evolution
equation 
\begin{equation}
	\partial _{t}u=\Delta _{p}u^{q}  \label{evoeq}
\end{equation}%
where $p>1$, $q>0$, $u=u(x,t)$ is an unknown function, and $%
\Delta _{p}$ is the $p$-Laplacian 
\begin{equation*}
	\Delta _{p}v=\func{div}\left( |\nabla v|^{p-2}\nabla v\right) .
\end{equation*}
For the physical meaning of (\ref{evoeq}) see \cite{grigor2023finite, leibenzon1945general, leibenson1945turbulent}. 

The equation (\ref{evoeq}) is referred to as a \textit{Leibenson} equation or a \emph{doubly non-linear parabolic} equation. In the case $q=1$, it becomes an evolutionary $p$-\textit{Laplace} equation $\partial _{t}u=\Delta_{p} u$,
and if in addition \(p=2\) then it amounts to the classical heat equation $\partial _{t}u=\Delta u$.

In the present paper, we are interested in the existence of solutions of the Cauchy-problem \begin{equation}
	\left\{ 
	\begin{array}{ll}
		\partial _{t}u=\Delta _{p}u^{q} &\text{in}~M\times (0, \infty), \\u(x, 0)=u_{0}(x)& \text{in}~M,
	\end{array}%
	\right.  \label{wholemandtvint}\end{equation}
where $u_{0}$ is a given function. 

We understand solutions of  (\ref{wholemandtvint}) in a certain weak sense (see Section \ref{secweak}
for the definition).
The first main result of the present paper is as follows  (cf. Theorem \ref{mainexthm}). 	

\begin{theorem}\label{mainthmint}
Assume that $u_{0}\in L^{1}(M)\cap L^{\infty}(M)$ is non-negative. Then there exists a unique non-negative bounded solution $u$ of the Cauchy problem (\ref{wholemandtvint}).
\end{theorem}

Let us emphasize that the only assumption that we impose on the manifold is geodesic completeness. 

Let us recall some previously known existence results for the Cauchy problem (\ref{wholemandtvint}).

In the special case when $M=\mathbb{R}^{n}$, existence results for the Cauchy problem (\ref{wholemandtvint}) were obtained for the whole range $q>0, p>1$, for example, in \cite{ishige1996existence, ivanov1997regularity, ladyzhenskaya1968linear, raviart1970resolution}.

In order to further discuss existence results on Riemannian manifolds, let us distinguish between the cases $q(p-1)>1$ (\textit{slow diffusion case}) and $q(p-1)<1$ (\textit{fast diffusion case}).

In the slow diffusion case, existence of solutions of (\ref{wholemandtvint}) was claimed by D. Andreucci and A. F. Tedeev in \cite{andreucci2015optimal} when $M$ has non-negative Ricci curvature. In fact, it was claimed therein under the weaker condition when $M$ satisfies a \textit{relative Faber-Krahn inequality} \cite{Buser, grigor, Saloff}. 
When also $p=2$, i.e. when (\ref{evoeq}) becomes the \textit{porous medium equation} with $q>1$, existence of solutions of (\ref{wholemandtvint}) was proved by N. De Ponti, G. Grillo, M. Muratori, C. Orrieri, F. Punzo and J. L. Vázquez in \cite{de2022wasserstein, grillo2016smoothing, grillo2018porous, grillo2017porous, vazquez2015fundamental} when $M$ is a \textit{Cartan-Hadamard manifold} and by G. Grillo, D. D. Monticelli and F. Punzo in \cite{grillo2025porous} when $M$ has non-negative Ricci curvature.

In the fast diffusion case, existence of solutions of (\ref{wholemandtvint}) for the porous medium equation ($q<1$) on Cartan-Hadamard manifolds was obtained in \cite{bonforte2008fast} by M. Bonforte, G. Grillo, and J. L. Vazquez and in \cite{grillo2021fast} by G. Grillo, M. Muratori, and F. Punzo.

Hence, the novelty of the present paper is that we prove existence without any further assumptions on the manifold besides the geodesic completeness. In particular, the existence is new in the case $p\ne 2$, $\frac{1}{p}\leq q<\frac{1}{p-1}$, (fast diffusion case) on any manifold besides $\mathbb{R}^{n}$.

Let us now discuss the approach of the proof of Theorem \ref{mainthmint}.
We first prove the existence of a solution of the mixed problem
\begin{equation}
	\left\{ 
	\begin{array}{ll}
		\partial _{t}u=\Delta_{p}u^{q} &\text{in}~B\times (0, T), \\ 
		u=0, & \text{on}~\partial B\times (0, T),\\u(x, 0)=u_{0}(x)& \text{in}~B,
	\end{array}%
	\right.  \label{dtvauxintlim}
\end{equation}
where $B$ is a precompact geodesic ball in $M$ and $T\in (0, \infty)$.
The method is inspired by \cite{bogelein2021existence} (see also \cite{ivanov1992existence}), where it was used to solve (\ref{dtvauxintlim}) when $B$ is a ball in $\mathbb{R}^{n}$ and when $p, q>1$.
For that, we consider, for $N>1$, the auxiliary problem \begin{equation}
	\left\{ 
	\begin{array}{ll}
		\partial _{t}u=\text{div}\left(A(u, \nabla u)\right) &\text{in}~B\times (0, T), \\ 
		u=\frac{1}{N}, & \text{on}~\partial B\times (0, T),\\u(x, 0)=u_{0}(x)+\frac{1}{N}& \text{in}~B,
	\end{array}%
	\right.  \label{dtvauxint}
\end{equation}%
where $A(u, \nabla u)=q^{p-1}\chi(u)^{(q-1)(p-1)}|\nabla u|^{p-2}\nabla u$ and $\chi(u)=\min\{N, \max\{u, \frac{1}{N}\}\}$  Using a modification of Galerkin's method (see also \cite{ladyzhenskaya1968linear}), we obtain solutions $u_{N, B, T}$ of (\ref{dtvauxint}). In the next step, we then prove that for large enough $N$, $\frac{1}{N}<u_{N, B, T}<C$, where $C$ does not depend on $N$, which implies that $A(u_{N, B, T}, \nabla u_{N, B, T})=|\nabla u_{N, B, T}^{q}|^{p-2}\nabla u_{N, B, T}^{q}.$ By some compactness argument, we then obtain for $N\to \infty$ a limit function $u_{B, T}$ of some subsequence of $u_{N, B, T}$, which is a solution of (\ref{dtvauxintlim}).
Proving suitable a priori norm estimates for solutions of (\ref{dtvauxintlim}) we are then able to send $B\to M$ and prove existence of a solution up to time $T$. By some uniqueness result, we then conclude the existence of solutions of (\ref{wholemandtvint}) for all $t\in (0, +\infty)$. 

The second part of the present paper deals with the property of \textit{finite propagation speed} of solutions of (\ref{wholemandtvint}). In \cite{grigor2023finite} it was proved that, under condition \begin{equation}\label{opticondint}q(p-1)>1,
\end{equation} solutions of (\ref{wholemandtvint}) have  finite propagation speed on an arbitrary Riemannian manifold (in the case \(q=1\) this was also proved in \cite{dekkers2005finite}). Under the additional restrictions 
\begin{equation}\label{sharprateintro1}p>2,\quad \frac{1}{p-1}<q\leq 1\end{equation} the estimate of the  rate of propagation for solutions of (\ref{wholemandtvint}) was improved in \cite{Grigoryan2024}.
Moreover, if in addition 
\begin{equation}\label{sharprateintro2}q< \frac{2}{p-1},\end{equation} it was shown in \cite{Grigoryan2024} that this estimate of the propagation is sharp for a class of spherically symmetric (model) manifolds (including $\mathbb{R}^n$).

The purpose of the second part of this paper is to obtain the improved propagation rate from \cite{Grigoryan2024} for a larger range of $p$ and $q$.
The second main result of the present paper is as follows  (cf. Theorem \ref{TFPSq}).

\begin{theorem}
	\label{mainthmintprop}         Assume that (\ref{opticondint}) holds. Let $u$ be a non-negative bounded solution to (\ref{wholemandtvint}). Let $\sigma $ be a real such that 
	\begin{equation}
		\sigma \geq 1\text{ and }\sigma >q(p-1)-1.  \label{si>>intro}
	\end{equation} If $u_{0}$ vanishes in a geodesic ball 
	$B_{0}$ in \(M\) of radius $R$ then 
	\begin{equation*}
		u=0\quad \text{\textnormal{in}}~\frac{1}{2}B_{0}\times \lbrack 0,t_{0}],
	\end{equation*}%
	where 
	\begin{equation}\label{tzintroop}
		t_{0}=\eta \mu(B_{0})^{\frac{q(p-1)-1}{\sigma}}R^{p}||u_{0}||_{L^{\sigma }(M)}^{-[q(p-1)-1]},
	\end{equation}%
	and the constant $\eta >0$ depends on the intrinsic geometry of $B_{0}$.
\end{theorem}
Hence, the solution $u$ has a finite propagation speed inside ball $B_{0}$, and
the rate of propagation is determined by \(t_{0}\) that depends on the intrinsic geometry of $B_{0}$ via
constant $\eta $.

The same result was obtained in \cite{Grigoryan2024} but under the stronger assumption (\ref{sharprateintro1}).

As in \cite{Grigoryan2024}, our Theorem \ref{mainthmintprop}  leads to a sharp propagation rate in a class of spherically symmetric manifolds provided $\sigma=1$ (cf. Corollary \ref{Corpropric}). By (\ref{si>>intro}), we can take $\sigma =1$ if $q(p-1)-1<1$ (which is equivalent to (\ref{sharprateintro2})) so that we obtain here a sharp propagation rate on these manifolds when $$\frac{1}{p-1}<q< \frac{2}{p-1}.$$ 

Let us discuss the differences in methods of the proof of finite propagation speed in \cite{Grigoryan2024} and the present paper and how they yield different ranges for $p$ and $q$. Even though, in both papers, the finite propagation speed follows from the same non-linear mean value inequality for subsolutions Lemma \ref{Tmeansing} and a modification of the classical De Giorgi iteration argument \cite{de1957sulla}, the ranges of $p$ and $q$ for which the mean value inequality holds are different. The mean value inequality in \cite{Grigoryan2024} is proved in the case (\ref{sharprateintro1}). This is because the proof uses the fact that, if $u$ is a non-negative subsolution of (\ref{wholemandtvint}), then the function 
\begin{equation}\label{utrun}
	(u^{a}-\theta)_{+}^{1/a}
\end{equation} 
is also a subsolution of (\ref{wholemandtvint}), provided $\theta\geq 0$ and 
\begin{equation}\label{condforaint}a:=\frac{q(p-1)-1}{p-2}\in (0, 1].\end{equation} 
In particular, the condition \(a\in(0,1]\) in (\ref{condforaint}) is satisfied provided (\ref{sharprateintro1}) holds. The second ingredient in this proof is a \textit{Caccioppoli-type inequality}, which says the following. Assume that $u$ is a subsolution of (\ref{wholemandtvint}) in $B\times I$ such that $u(\cdot, 0)=0$ in a ball $B$. Let $\eta \left( x,t\right) $ be a locally Lipschitz non-negative bounded function. Then for any $\lambda$ large enough, $\sigma=\lambda+q(p-1)-1$, 
$$\left[ \int_{B }u^{\lambda }\eta ^{p}\right] _{0}^{t}+c_{1}%
\int_{B\times [0, t]}\left\vert \nabla \left( u^{\sigma/p }\eta \right) \right\vert
^{p}\leq c_{2}\int_{B\times [0, t]}u^{\sigma }\left\vert \nabla \eta \right\vert ^{p}.$$ Then the aforementioned property allows to apply this inequality to the subsolutions $$u_{k}=\left( u^{a}-\left( 1-2^{-k}\right) \theta \right) _{+}^{1/a},\quad k\geq 0, $$ for some fixed $\theta>0$, where $a$ is given by (\ref{condforaint}).

However, in the present paper we prove the mean value inequality for the whole range $p>1, q>0$. This is because we use in the proof instead a Caccioppoli type inequality of the form $$\left[ \int_{B }u_{k+1}^{\lambda }\eta ^{p}\right] _{0}^{t}+A^{k}%
\int_{B\times [0, t]}\left\vert \nabla \left( u_{k+1}^{\sigma/p }\eta \right) \right\vert
^{p}\leq B^{k}\int_{B\times [0, t]}u_{k}^{\sigma }\left\vert \nabla \eta \right\vert ^{p}.$$ where $$u_{k}=\left( u-\left( 1-2^{-k}\right) \theta \right) _{+},\quad k\geq 0,$$ and $A, B$ are positive constants (cf. Lemma \ref{Lem1}).

The structure of the paper is as follows.

In Section \ref{secweak}, we
define the notion of a weak solution of the Leibenson equation (\ref{evoeq}) and prove a \textit{comparison principle} following a method from \cite{anttila2026uniqueness, otto1996l1}, which later yields the uniqueness.

In Section \ref{secaux} we prove in Lemma \ref{lemexist} the existence of solutions of auxiliary problem (\ref{dtvauxint}).

In Section \ref{secmvi} we prove the central technical
result for the second part of this paper $-$ the mean value inequality for subsolutions (Lemma \ref{Tmeansing}).

In Section \ref{secglobex} we prove our main existence Theorem \ref{mainthmint}. 

In Section \ref{secfps} we prove our
main results about finite propagation speed.

For other qualitative and quantitative properties of solutions of (\ref{evoeq}) on Riemannian manifolds see \cite{grigor2025upper, meglioli2025global, surig2024finite, surig2024sharp, surig2025gradient} and from a probabilistic point of view \cite{barbu2025leibenson}.

\section{Weak solutions}\label{secweak}
\subsection{Definition and basic properties}

Let $M$ be a Riemannian manifold and $\Omega$ be an open subset of $M$. We consider in what follows the following Cauchy problem:%
\begin{equation}
	\left\{ 
	\begin{array}{ll}
		\partial _{t}u=\Delta _{p}u^{q} &\text{in}~\Omega\times (0, T), \\ 
		u=0, & \text{on}~\partial \Omega\times (0, T),\\u(x, 0)=u_{0}(x)& \text{in}~\Omega,
	\end{array}%
	\right.  \label{dtv}
\end{equation}%
where $u_{0}$ is a given function, in a certain weak sense as explained below.

We assume throughout that 
\begin{equation*}
	p>1\ \ \text{and}\ \ \ q>0.
\end{equation*}

Let $\mu $ denote the Riemannian measure on $M$. For simplicity of notation,
we frequently omit in integrations the notation of measure. All integration
in $M$ is done with respect to $d\mu $, and in $M\times \mathbb{R}$ -- with
respect to $d\mu dt$, unless otherwise specified.

\begin{definition}\label{defweaksolu}
	\normalfont
	We say that a function $u=u(x, t)$ is a \textit{weak subsolution} of (\ref{dtv}) in $\Omega\times (0, T)$, where $0<T< \infty$, if 
	\begin{equation}  \label{defvonsoluq}
		u\in C\left([0, T]; L^{q+1}(\Omega)\right) \text{\ and\ }
		u^{q}\in L^{p}\left((0, T); D_{0}^{1, p}(\Omega)\right)
	\end{equation}
	and for all non-negative \textit{test functions} $\psi$ such that
	\begin{equation}  \label{defvontestsoluq}
		\psi\in W^{1, 1+\frac{1}{q}}\left((0, T);L^{1+\frac{1}{q}}(\Omega)\right)
		\cap L^{p}\left((0, T); D_{0}^{1, p}(\Omega)\right)
		,
	\end{equation}
	and for all $t_{1}, t_{2}\in [0, T]$ with $t_{1}<t_{2}$, we have 
	\begin{equation}  \label{defvonweaksolq}
		\left[\int_{\Omega}{u\psi}\right]_{t_{1}}^{t_{2}}+\int_{t_{1}}^{t_{2}}{%
			\int_{\Omega}{-u\partial_{t}\psi+|\nabla u^{q}|^{p-2}\langle\nabla u^{q},
				\nabla \psi\rangle}}\leq 0.
	\end{equation}

Here, $D_{0}^{1, p}(\Omega)$ denotes the completion of $C_{0}^{\infty}(\Omega)$ with respect to the seminorm $||u||_{D_{0}^{1, p}(\Omega)}=||\nabla u||_{L^{p}(\Omega)}$.

The condition $u(x, 0)=u_{0}(x)$ is understood in the sense of $L^{q+1}(\Omega)$, meaning that, $$\lim_{t\to 0}\int_{\Omega}|u(\cdot, t)-u_{0}|^{q+1}=0.$$
The notion of a \textit{weak supersolution} and \textit{weak solution} is defined analogously.
\end{definition}

Note that a priori we do not assume that $u$ is non-negative. In the case when $u$ is not non-negative, we define $u^{q}:=|u|^{q-1}u$. However, since we assume that $u_{0}$ is non-negative, this will imply that $u$ is non-negative (cf. Lemma \ref{thmcomp}).

From now on let $0<T<\infty$.

\begin{definition}
	\normalfont
	Let $u=u(x, t)\in L^{1}(\Omega\times (0, T))$ and $u(\cdot,
	0)=u_{0}$. Then we define, for $h\in (0, T)$, 
	\begin{equation}\label{defsteklov}
		u_{h}(\cdot, t)=\frac{1}{h}\int_{0}^{t}{e^{(s-t)/h}u(\cdot, s)ds}
	\end{equation}
	and 
	\begin{equation*}u_{\overline{h}}(\cdot, t)=\frac{1}{h}\int_{t}^{T}{e^{(t-s)/h}u(\cdot,
			s)ds}.\end{equation*}
We also define $u_{h_{0}}(\cdot, t)=e^{-t/h}u_{0}+u_{h}(\cdot, t)$.
\end{definition}

The properties of $u_{h_{0}}$ and $u_{h}$ stated in the following Lemma \ref{convstek} are proved in
Lemma 2.2 in \cite{kinnunen2006pointwise} and in Lemma B.1 and Lemma B.2 in 
\cite{bogelein2013parabolic}.

\begin{lemma}\label{propstek}
	\label{convstek} 	
If $u\in L^{p}\left((0, T); W^{1, p}(\Omega)\right)$ then $u_{h}\in L^{p}\left((0, T); W^{1, p}(\Omega)\right)$ and $u_{h}\to u$ in $L^{p}\left((0, T); W^{1, p}(\Omega)\right)$ as $h\to 0$. A similar result also holds in the space $L^{p}\left((0, T); W_{0}^{1, p}(\Omega)\right)$.
	
	If $u\in L^{p}(\Omega\times (0, T))$, then 
	
	\begin{equation}  \label{propmol}
		\partial_{t}u_{h}=\frac{1}{h}(u-u_{h})\qquad \textnormal{and}\qquad \partial_{t}u_{\overline{h}}=\frac{1}{h}(u_{\overline{h}}-u).
	\end{equation}
If $u\in L^{p}\left((0, T); L^{p}(\Omega)\right)$, then $u_{h}, u_{\overline{h}}\in C([0, T], L^{p}(\Omega)).$ Finally, if $u\in C([0, T], L^{p}(\Omega))$, we also have $u_{h}(\cdot, t)\to u(\cdot, t)$ in $L^{p}(\Omega)$ as $h\to 0$ and the same holds for $u_{\overline{h}}$.
\end{lemma}

\begin{lemma}\label{continuitylemma}
Let $Q=\Omega\times (0, T)$. Let $u=u(x, t)\in L^{q+1}(Q)$ be a non-negative function such that $u^{q}\in L^{p}\left((0, T); W_{0}^{1, p}(\Omega)\right)$ satisfying \begin{equation}\label{condforcont}\int_{Q}{-u\partial_{t}\psi+|\nabla u^{q}|^{p-2}\langle\nabla u^{q},\nabla \psi\rangle}= 0,\end{equation} where $\psi$ is as in (\ref{defvontestsoluq}) with $\psi(0)=\psi(T)=0$. Then $u$ has a representative in $C([0, T]; L^{q+1}(\Omega))$.
\end{lemma}

\begin{proof}
Let $w\in L^{q+1}(Q)$ satisfy $w^{q}\in L^{p}\left((0, T); W_{0}^{1,p}(\Omega)\right)$ and $\partial_{t}w^{q}\in L^{\frac{q+1}{q}}(Q)$. We want to test in (\ref{condforcont}) with $\psi=\zeta(w^{q}-(u^{q})_{h})$, where $w^{q}=(u^{q})_{\overline{h}}$ and $\zeta(t)=\eta(t)\varphi(t)$, where $\eta(t)=\frac{1}{\varepsilon}(t-\tau)$ if $t\in [\tau, \tau+\varepsilon]$ and $\varphi(0)=\varphi(T)=0$ and $|\varphi^{\prime}(t)|\leq \frac{C}{T}$. We have \begin{align*}
\int_{Q}u\partial_{t}\psi&=\int_{Q}u\zeta^{\prime}(w^{q}-(u^{q})_{h})+u\zeta(\partial_{t}w^{q}-\partial_{t}(u^{q})_{h})\\&=\int_{Q}u\zeta^{\prime}(w^{q}-(u^{q})_{h})+u\zeta \partial_{t}w^{q}-\zeta (u^{q})_{h}^{1/q}\partial_{t}(u^{q})_{h}+\zeta((u^{q})_{h}^{1/q}-u)\partial_{t}(u^{q})_{h}\\&=\int_{Q}u\zeta^{\prime}(w^{q}-(u^{q})_{h})+u\zeta \partial_{t}w^{q}-\zeta \frac{q}{q+1}\partial_{t}(u^{q})_{h}^{\frac{q}{q+1}}-\frac{1}{h}\zeta((u^{q})_{h}^{1/q}-u)((u^{q})_{h}-u^{q})\\&\leq \int_{Q}u\zeta^{\prime}(w^{q}-(u^{q})_{h})+u\zeta \partial_{t}w^{q}-\zeta \frac{q}{q+1}\partial_{t}(u^{q})_{h}^{\frac{q}{q+1}},
\end{align*} where we used (\ref{propmol}) in the penultimate step. Further, we have \begin{align*}
\int_{Q}u\partial_{t}\psi\leq \int_{Q}u\zeta^{\prime}(w^{q}-(u^{q})_{h})+u\zeta \partial_{t}w^{q}+\zeta^{\prime}\frac{q}{q+1}(u^{q})_{h}^{\frac{q}{q+1}}.
\end{align*}
Hence, we deduce \begin{align*}
\lim_{h\to 0}\int_{Q}u\partial_{t}\psi&\leq \int_{Q}u\zeta^{\prime}(w^{q}-u^{q})+u\zeta \partial_{t}w^{q}+\zeta^{\prime}\frac{q}{q+1}u^{q+1}\\&=\int_{Q}\zeta(u-w)\partial_{t}w^{q}+\zeta\frac{q}{q+1}\partial_{t}w^{q+1}+\zeta^{\prime}\left(\frac{q}{q+1}u^{q+1}+u(w^{q}-u^{q})\right)\\&=\int_{Q}\zeta(u-w)\partial_{t}w^{q}-\zeta\frac{q}{q+1}w^{q+1}+\zeta^{\prime}b(u, w),
\end{align*}
where $b(u, w)=\frac{q}{q+1}(w^{q+1}-u^{q+1})-u(w^{q}-u^{q})$.
On the other hand, we get from Lemma \ref{propstek}, $$\lim_{h\to 0}\int_{Q}{|\nabla u^{q}|^{p-2}\langle\nabla u^{q},\nabla (w^{q}-(u^{q})_{h})\rangle\zeta}=\int_{Q}{|\nabla u^{q}|^{p-2}\langle\nabla u^{q},\nabla (w^{q}-u^{q})\rangle\zeta}.$$
Hence, it follows that \begin{align}\label{equalityforcont}\int_{Q}\zeta^{\prime}b(u, w)\leq\int_{Q}\zeta(\partial_{t}w^{q}(u-w)+|\nabla u^{q}|^{p-2}\langle\nabla u^{q},\nabla u^{q}-\nabla w^{q}\rangle).\end{align}
For later usage, let us also test with $\psi=\zeta(w^{q}-(u^{q})_{\overline{h}})$ in (\ref{defvonweaksolq}) so that by using similar arguments, we obtain the opposite inequality, which yields that \begin{align}\label{laterusagefromcont}\int_{Q}\zeta^{\prime}b(u, w)=\int_{Q}\zeta(\partial_{t}w(u-w)+|\nabla u^{q}|^{p-2}\langle\nabla u^{q},\nabla u^{q}-\nabla w\rangle).\end{align}
Applying inequality (\ref{equalityforcont})  we get using (\ref{propmol}) \begin{align*}\frac{1}{\varepsilon}\int_{\tau}^{\tau+\varepsilon}\int_{\Omega}b(u, (u^{q})_{\overline{h}}^{\frac{1}{q}})&\leq \frac{C}{T}\int_{Q}b(u, (u^{q})_{\overline{h}}^{\frac{1}{q}})\\&\quad+\int_{Q}\frac{1}{h}((u^{q})_{\overline{h}}-u^{q})(u-(u^{q})_{\overline{h}}^{\frac{1}{q}})+|\nabla u^{q}|^{p-2}|\langle\nabla u^{q},\nabla u^{q}-\nabla (u^{q})_{\overline{h}}\rangle|\\&\leq\frac{C}{T}\int_{Q}b(u, (u^{q})_{\overline{h}}^{\frac{1}{q}})+|\nabla u^{q}|^{p-1}|\nabla u^{q}-\nabla (u^{q})_{\overline{h}}|.
\end{align*}
Let us now use the estimates $$c|x-y|^{q+1}\leq b(x, y)\leq C|x^{q}-y^{q}|^{\frac{q+1}{q}},$$ which holds for all $x, y\geq 0$, where $c$ and $C$ only depend on $q$ (see Lemmas 2.2 \& 2.3 in \cite{bogelein2018higher}). Sending $\varepsilon \to 0$ it follows that $$\textnormal{ess}\sup_{0<\tau<T}\int_{\Omega}\left|u- (u^{q})_{\overline{h}}^{\frac{1}{q}}\right|^{q+1}(\cdot, \tau)\leq \frac{C}{T}\int_{Q}|u^{q}- (u^{q})_{\overline{h}}|^{\frac{q+1}{q}}+|\nabla u^{q}|^{p-1}|\nabla u^{q}-\nabla (u^{q})_{\overline{h}}|.$$ Hence, we get \begin{equation}\label{repest}\lim_{h\to 0}\textnormal{ess}\sup_{0<\tau<T}\int_{\Omega}\left|u- (u^{q})_{\overline{h}}^{\frac{1}{q}}\right|^{q+1}(\cdot, \tau)=0.\end{equation} Note that the function $w=(u^{q})_{\overline{h}}^{\frac{1}{q}}:[0, T]\to L^{q+1}(\Omega)$ is continuous as $$|w(x, s)-w(x, t)|^{q+1}\leq |w(x, s)^{q}-w(x, t)^{q}|^{\frac{q+1}{q}}=|(u^{q})_{\overline{h}}(x, s)-(u^{q})_{\overline{h}}(x, t)|^{\frac{q+1}{q}}$$ and $(u^{q})_{\overline{h}}:[0, T]\to L^{q+1}(\Omega)$ is continuous by Lemma \ref{propstek}. From (\ref{repest}) we therefore conclude that $u$ has a continuous representative.
\end{proof}

\subsection{Comparison principle}\label{uniqe}

	\begin{lemma}\label{thmcomp} Let $u$ and $v$ be two weak solutions of (\ref{evoeq}) in $M\times (0, T)$ with inital functions $u_{0}, v_{0}$ respectively. Assume that
	\begin{equation}\label{initial}
		u_{0} \leq v_{0}\quad \textnormal{a.e. in}~M. 
	\end{equation}
	Then $u \leq v$ a.e. in $M\times (0, T)$.
\end{lemma}

The proof of Lemma \ref{thmcomp} follows from an appropriate choice of test function in the weak formulation (\ref{defvonweaksolq}). For that, let $\eta \in C^\infty(\mathbb{R})$ be a convex function satisfying
\[
\eta(z) :=
\begin{cases}
	0 & \text{when } z \le 0, \\
	z - \frac{1}{2} & \text{when } z \ge 1.
\end{cases}
\]
Note that the derivative $\eta'$ approximates the Heaviside function and converges to the sign-function. For $\sigma > 0$, let us also define 
\[
\eta_{\sigma,+}(z) := \sigma \eta\!\left(\frac{z}{\sigma}\right), \quad
\eta_{\sigma,-}(z) := \sigma \eta\!\left(-\frac{z}{\sigma}\right),
\]
and
\begin{equation}
	\varphi_{\sigma,\pm}(z, w) := \int_w^z \eta'_{\sigma,\pm}(s^{q} - w^{q}) \, ds, 
\end{equation}
where $\eta'_{\sigma,\pm}$ is the derivative of $\eta_{\sigma,\pm}$.

\begin{lemma}
	For all $z_1, z_2, w \in \mathbb{R}$ and $\sigma > 0$ we have
	\begin{equation}\label{keylemma}
		\varphi_{\sigma,\pm}(z_1, w) - \varphi_{\sigma,\pm}(z_2, w) \ge \eta'_{\sigma,\pm}(z_2^{q} - w^{q}) (z_1 - z_2).\end{equation}
\end{lemma}

\begin{proof}
	Since $\eta$ is a smooth, convex function, the same is true for $\eta_{\sigma,\pm}$ and thus $\eta_{\sigma,\pm}^{\prime}$ is non-decreasing. Since $s\mapsto s^{q}$ is also non-decreasing, we see that $s\mapsto \eta_{\sigma,\pm}^{\prime}(s^{q})$ is non-decreasing.
	Therefore, we have in the case $z_1 \ge z_2$
	\begin{align*}
		\varphi_{\sigma,+}(z_1, w) - \varphi_{\sigma,+}(z_2, w)
		= \int_{z_2}^{z_1} \eta'_{\sigma,+}(s^{q} - w^{q}) ds
		\ge \eta'_{\sigma,+}(z_2^{q} - w^{q}) (z_1 - z_2).
	\end{align*}
	If $z_1 < z_2$, the same estimate holds as
	$$\varphi_{\sigma,+}(z_1, w) - \varphi_{\sigma,+}(z_2, w)
	= -\int_{z_1}^{z_2} \eta'_{\sigma,+}(s^{q} - w^{q})\, ds
	\ge \eta'_{\sigma,+}(z_2^{q} - w^{q}) (z_1 - z_2).$$
	The proof for $\varphi_{\sigma,-}$ follows by similar arguments so that we conclude (\ref{keylemma}).
\end{proof} 

Let us also fix a function $\phi \in C^\infty_0(\mathbb{R})$ satisfying
$$
\begin{cases}
	\phi(z) \ge 0 \text{ for all } z \in \mathbb{R}, \\
	\phi(z) = 0 \text{ for all } |z| \ge 1, \\
	\int_{\mathbb{R}} \phi\, dt = 1.
\end{cases}$$
For $\epsilon > 0$ and a non-negative function $\gamma \in C^\infty_0(0, T)$ we define $\gamma_\epsilon \in C^\infty_0(\mathbb{R}^2)$ by
\begin{equation}\label{defgamma}
	\gamma_\epsilon(t, s) := \frac{1}{\epsilon} \phi\!\left(\frac{t - s}{\epsilon}\right)\gamma\!\left(\frac{t + s}{2}\right).
\end{equation}
Later we will send $\epsilon$ to $0$.

In the following we sometimes omit for simplicity the spatial variable $x$ and write $u(t) = u(x, t)$.

Fix some precompact open subset $\Omega$ of $M$.

\begin{lemma}
	Let $u$ and $v$ be two weak solutions of (\ref{dtv}) with inital functions $u_{0}, v_{0}$ respectively. Let $\omega\in C_{0}^{\infty}(\Omega)$ be non-negative. Then, for all $\epsilon>0$ small enough, all $\sigma > 0$, and almost every $s \in (0, T)$, we have
	\begin{equation}
		\begin{aligned}
			\int_0^T \int_\Omega \Big(
			&-\varphi_{\sigma,+}(u(t), v(s)) \partial_t \gamma_\epsilon(t, s)\omega \\
			&+ \gamma_\epsilon(t, s)\eta''_{\sigma,+}(u^{q}(t)-v^{q}(s))
			\langle |\nabla u^{q}(t)|^{p-2} \nabla u^{q}(t), \nabla u^{q}(t)-\nabla v^{q}(s)\rangle\omega\\
			&\label{test+}+\gamma_\epsilon(t, s)\eta'_{\sigma,+}(u^{q}(t)-v^{q}(s))
			\langle |\nabla u^{q}(t)|^{p-2} \nabla u^{q}(t), \nabla \omega \rangle
			\Big)\, dxdt\leq 0.
		\end{aligned}
	\end{equation}
	Similarly for small enough $\epsilon$, all $\sigma > 0$, and almost every $t \in (0, T)$, it holds that
	\begin{equation}
		\begin{aligned}
			\int_0^T \int_\Omega \Big(
			&-\varphi_{\sigma,-}(v(s), u(t)) \partial_s \gamma_\epsilon(t, s)\omega \\
			&- \gamma_\epsilon(t, s)\eta''_{\sigma,+}(u^{q}(t)-v^{q}(s))
			\langle |\nabla v^{q}(s)|^{p-2} \nabla v^{q}(s), \nabla u^{q}(t)-\nabla v^{q}(s)\rangle\omega\\
			&\label{test-}+\gamma_\epsilon(t, s)\eta'_{\sigma,+}(u^{q}(t)-v^{q}(s))
			\langle |\nabla v^{q}(t)|^{p-2} \nabla v^{q}(t), \nabla \omega \rangle
			\Big)\, dxds\leq 0.
		\end{aligned}
	\end{equation}
\end{lemma} 

\begin{proof}
	Note that for small enough $\epsilon > 0$, only depending on the support of $\gamma$, it holds that
	\begin{equation}\label{compsupp}
		t \mapsto \gamma_\epsilon(t, s) \in C^\infty_0(0, T)
	\end{equation}
	for all $s > 0$.
	Also observe that $$
	\psi(x, t) := \eta'_{\sigma,+}(u^{q}(x, t) - v^{q}(x, s)) \gamma_\epsilon(t, s)\omega(x)\in L^p(0, T; W_{0}^{1,p}(\Omega))
	$$
	for almost every $s \in (0, T)$.
	
	For notational convenience, we extend $\psi(\cdot, t) = 0$ for $t > T$ and $t < 0$. Let us consider the Steklov averages \begin{equation}\label{defclassicstekl}
	\psi^h(x, t) := \frac{1}{h} \int_t^{t+h} \psi(x, \tau)\, d\tau.\end{equation}
	Then it follows from \cite{bogelein2013parabolic} that $\psi^h \in L^p(0, T; W^{1,p}_0(\Omega))$ for all small enough $h > 0$ and that $\psi^h$ has the time derivative
	$$
	\partial_t \psi^h(t) = \frac{\psi(t+h) - \psi(t)}{h} \in L^\infty(M\times (0, T)).
	$$
	Hence, we can use $\psi^h$ as a test function in (\ref{defvonweaksolq}) and obtain
	$$
	\int_0^T \int_\Omega u(t) \partial_t \psi^h(t)\, dxdt
	=
	\int_0^T \int_\Omega \langle |\nabla u^{q}(t)|^{p-2} \nabla u^{q}(t), \nabla \psi^h(t) \rangle\, dxdt,
	$$ where we used that $\psi_{h}(0)=\psi_{h}(T)=0$ for small enough $h>0$ by (\ref{compsupp}).
	
	By a change of variables and the notational conventions concerning $\psi$ we obtain
	\begin{align*}
		\int_0^T \int_\Omega u(t) \partial_t \psi^h(t)\, dxdt
		&=
		\int_0^T \int_\Omega u(t) \frac{\psi(t+h) - \psi(t)}{h}\, dxdt
		\\&=
		\int_0^T \int_\Omega
		- \frac{u(t) - u(t-h)}{h} \psi(t)\, dxdt.	
	\end{align*}
	By the definition of $\psi$ and using (\ref{keylemma}), we further deduce
	\begin{align*}
		\int_0^T \int_\Omega \frac{u(t) - u(t - h)}{h} \psi(t)\, dxdt
		&=
		\int_0^T \int_\Omega
		\frac{u(t) - u(t - h)}{h}
		\eta'_{\sigma,+}(u^{q}(t)-v^{q}(s)) \gamma_\epsilon(t, s)\omega\, dxdt
		\\&\ge
		\int_0^T \int_\Omega
		\frac{\varphi_{\sigma,+}(u(t), v(s)) - \varphi_{\sigma,+}(u(t - h), v(s))}{h}
		\gamma_\epsilon(t, s)\omega\, dxdt
		\\&=
		\int_0^T \int_\Omega
		- \varphi_{\sigma,+}(u(t), v(s))
		\frac{\gamma_\epsilon(t + h, s) - \gamma_\epsilon(t, s)}{h}\omega
		\, dxdt,
	\end{align*}
	where we used in the last equality that $v(s)$ does not depend on $t$. Combining this with the computations above we get
	$$
	\int_0^T \int_\Omega
	- \varphi_{\sigma,+}(u(t), v(s))
	\frac{\gamma_\epsilon(t + h, s) - \gamma_\epsilon(t, s)}{h}\omega
	\, dxdt
	\le
	- \int_0^T \int_\Omega
	\langle |\nabla u^{q}(t)|^{p-2} \nabla u^{q}(t), \nabla \psi^h(t) \rangle\, dxdt.
	$$
	Note that \begin{align*}
		\nabla \psi(t)=\gamma_\epsilon(t, s)\eta''_{\sigma,+}(u^{q}(t)-v^{q}(s))(\nabla u^{q}(t) - \nabla v^{q}(s))\omega+\gamma_\epsilon(t, s)\eta'_{\sigma,+}(u^{q}(t)-v^{q}(s))\nabla \omega.
	\end{align*}
	By letting $h \to 0^+$ and using the convergence properties of the Steklov averages \cite{bogelein2013parabolic}, it follows that
	\begin{align*}
		0 &\ge \int_0^T \int_\Omega
		- \varphi_{\sigma,+}(u(t), v(s)) \partial_t \gamma_\epsilon(t, s)\omega
		+ \langle |\nabla u^{q}(t)|^{p-2} \nabla u^{q}(t), \nabla \psi(t) \rangle
		\, dxdt
		\\&=
		\int_0^T \int_\Omega
		- \varphi_{\sigma,+}(u(t), v(s)) \partial_t \gamma_\epsilon(t, s)\omega
		\\&\qquad+ \gamma_\epsilon(t, s)\eta''_{\sigma,+}(u^{q}(t)-v^{q}(s))
		\langle |\nabla u^{q}(t)|^{p-2} \nabla u^{q}(t), \nabla u^{q}(t) - \nabla v^{q}(s) \rangle\omega\,\\&\qquad+\gamma_\epsilon(t, s)\eta'_{\sigma,+}(u^{q}(t)-v^{q}(s))
		\langle |\nabla u^{q}(t)|^{p-2} \nabla u^{q}(t), \nabla \omega \rangle  dxdt,
	\end{align*}
	which proves (\ref{test+}).
	
	The proof of (\ref{test-}), follows from similar arguments by considering the test function
	$$\psi(x, s) = \eta_{\sigma,-}^{\prime}(v^{q}(s) - u^{q}(t)) \gamma_\epsilon(t, s)\omega\in L^p(0, T; W_{0}^{1,p}(\Omega))$$ for small enough $\epsilon > 0$, all $\sigma > 0$ and almost every $t \in (0, T)$.
\end{proof} 

\begin{lemma}\label{plapmon}
	Let $p>1$. Then, for any $\xi, \eta \in T_{x}M$, \begin{equation}\label{monotoplap}\langle |\xi|^{p-2}\xi- |\eta|^{p-2}\eta, \xi-\eta\rangle\geq c\left\{ 
		\begin{array}{ll}
			|\xi-\eta|^{p} &\text{if}~p\geq2, \\ 
			\left(|\xi|^{2}+|\eta|^{2}\right)^{\frac{p-2}{2}}|\xi-\eta|^{2}, & \text{if}~p<2,
		\end{array}%
		\right.\end{equation}
	where $c$ is a positive constant depending on $p$ and we assume in the latter case that $\xi$ or $\eta$ are non-zero.
\end{lemma}

In the next lemma we obtain estimates by sending $\sigma \to 0+$ in (\ref{test+}) and (\ref{test-}).

\begin{lemma}
	Let $u$ and $v$ be two weak solutions of (\ref{dtv}) with inital functions $u_{0}, v_{0}$ respectively. Then, for all $\epsilon>0$ small enough, 
	\begin{equation}\label{addingup}
		\int_0^T \int_0^T \int_\Omega
		(u(t) - v(s))_+
		(\partial_t + \partial_s)\gamma_\epsilon(t, s)\omega+\gamma_\epsilon(t, s)(|\nabla u^{q}(t)|^{p-1}+ |\nabla v^{q}(t)|^{p-1}) |\nabla \omega|\, dx\,dt\,ds \ge 0.
	\end{equation}
\end{lemma}

\begin{proof}
	Integrating (\ref{test+}) over $s \in (0, T)$ and (\ref{test-}) over $t \in (0, T)$, and adding these two integrals, it follows that
	\begin{align*}
		\int_0^T \int_0^T \int_\Omega
		&- \varphi_{\sigma,+}(u(t), v(s)) \partial_t \gamma_\epsilon(t, s)\omega
		- \varphi_{\sigma,-}(v(t), u(s)) \partial_s \gamma_\epsilon(t, s)\omega
		\\&+ \gamma_\epsilon(t, s)\eta''_{\sigma,+}(u(t)-v(s)) \mathcal{A}(\nabla u^{q}(t), \nabla v^{q}(s))\\&-\gamma_\epsilon(t, s)\eta'_{\sigma,+}(u^{q}(t)-v^{q}(s))(|\nabla u^{q}(t)|^{p-1}+ |\nabla v^{q}(t)|^{p-1}) |\nabla \omega|
		\, dxdtds\leq 0,
	\end{align*}
	where we used the notation
	$$\mathcal{A}(\nabla u^{q}, \nabla v^{q})
	=
	\left\langle
	|\nabla u^{q}|^{p-2}\nabla u^{q} - |\nabla v^{q}|^{p-2}\nabla v^{q},
	\nabla u^{q} - \nabla v^{q}
	\right\rangle
	.$$
	Since $\eta''_{\sigma,+} \ge 0$, $\gamma_\epsilon \ge 0$ and (\ref{monotoplap}), it follows that the term on the second row is non-negative.
	Hence, \begin{align}\nonumber
		\int_0^T \int_0^T \int_\Omega
		&\varphi_{\sigma,+}(u(t), v(s)) \partial_t \gamma_\epsilon(t, s)\omega
		+ \varphi_{\sigma,-}(v(t), u(s)) \partial_s \gamma_\epsilon(t, s)\omega\\&\label{beforesig}+\gamma_\epsilon(t, s)\eta'_{\sigma,+}(u^{q}(t)-v^{q}(s))(|\nabla u^{q}(t)|^{p-1}+ |\nabla v^{q}(t)|^{p-1}) |\nabla \omega|
		\, dxdtds \ge 0.
	\end{align}
	
	We now want to let $\sigma \to 0+$ in (\ref{beforesig}). First, note that when $\sigma \to 0+$,
	\[
	\eta'_{\sigma,+}(z) \to \operatorname{sgn}(z) :=
	\begin{cases}
		0, & z \le 0,\\
		1, & z > 0,
	\end{cases}
	\]
	uniformly in $\mathbb{R} \setminus [-a,a]$ where $a > 0$ is fixed. Hence, as $\sigma \to 0+$, we have
	$\varphi_{\sigma,+}(z,w) \to (z - w)_+ $ pointwise in $\mathbb{R}^2$.
	Clearly, we also have
	$|\varphi_{\sigma,+}(z,w)| \le (z - w)_+$ for all $z,w \in \mathbb{R}$ and  $\sigma > 0.$ Analogously it also holds that
	$\varphi_{\sigma,-}(z,w) \to (w - z)_+$ and $|\varphi_{\sigma,-}(z,w)| \le (z - w)_+$.
	Finally, (\ref{addingup}) follows from the dominated convergence theorem by sending $\sigma \to 0+$ in (\ref{beforesig}), which was to be proved
\end{proof}	

In the next lemma we let $\epsilon \to 0+$ in (\ref{addingup}).

\begin{lemma}
	Let $u$ and $v$ be two weak solutions of (\ref{dtv}) with inital functions $u_{0}, v_{0}$ respectively. Then
	\begin{equation}\label{epsi0}
		\int_0^T \int_\Omega (u(t) - v(t))_+ \gamma'(t)\omega+(|\nabla u^{q}(t)|^{p-1}+ |\nabla v^{q}(t)|^{p-1}) |\nabla \omega|\, dxdt \ge 0.
	\end{equation}
\end{lemma} 

\begin{proof}
	Let $\epsilon > 0$ be small enough so that (\ref{addingup}) holds and that $\gamma_\epsilon \in C^\infty_0((0, T)\times(0, T))$. From (\ref{defgamma}) we get
	$$
	(\partial_t + \partial_s)\gamma_\epsilon(t, s)
	=
	\frac{1}{\epsilon}
	\phi\!\left(\frac{t - s}{\epsilon}\right)
	\gamma'\!\left(\frac{t + s}{2}\right).
	$$
	Therefore, 
	\begin{align*}
		\int_0^T \int_0^T \int_\Omega
		&(u(t) - v(s))_+
		\frac{1}{\epsilon}
		\phi\!\left(\frac{t - s}{\epsilon}\right)
		\gamma'\!\left(\frac{t + s}{2}\right)\omega\\&+\gamma_\epsilon(t, s)(|\nabla u^{q}(t)|^{p-1}+ |\nabla v^{q}(t)|^{p-1}) |\nabla \omega|
		\, dxdtds \ge 0.
	\end{align*}
	
	By the choice of $\phi$ and by letting $\epsilon \to 0+$ we conclude
	$$
	\int_0^T \int_\Omega (u(t) - v(t))_+ \gamma'(t)\omega+\gamma(t)(|\nabla u^{q}(t)|^{p-1}+ |\nabla v^{q}(t)|^{p-1}) |\nabla \omega|\, dxdt \ge 0,
	$$
	which finishes the proof.
\end{proof} 

\begin{proof}[Proof of Lemma \ref{thmcomp}] Let us choose $\gamma$ in (\ref{epsi0}) as a cut-off function, approximating the characteristic function of the interval $[t_1, t_2]$, where $0 < t_1 < t_2 < T$. Hence, we obtain in the limit and using that $u\omega$ and $v\omega\in C([0, T], L^{1}(\Omega))$,
	\begin{align*}
		\int_\Omega (u(t_2) - v(t_2))_+\omega\, dx+\int_{t_{1}}^{t_{2}} \int_\Omega (|\nabla u^{q}(t)|^{p-1}+ |\nabla v^{q}(t)|^{p-1}) |\nabla \omega|\, dxdt
		\le
		\int_\Omega (u(t_1) - v(t_1))_+\omega\, dx.
	\end{align*}
	For $t_{1}\to 0$ we get
	\begin{align*}
		\int_\Omega (u(t_2) - v(t_2))_+\omega\, dx+\int_{0}^{t_{2}} \int_\Omega (|\nabla u^{q}(t)|^{p-1}+ |\nabla v^{q}(t)|^{p-1}) |\nabla \omega|\, dxdt
		\le
		\int_\Omega (u_0 - v_0)_+\omega\, dx.
	\end{align*}
	Finally, choosing $\omega$ such that $\omega\to 1$ and $|\nabla \omega|\to 0$ in $\Omega$, we deduce by (\ref{initial}) for any $t\in (0, T)$, 
	\begin{equation}\label{forcomparison}\int_\Omega (u(t) - v(t))_+\, dx \leq \int_\Omega (u_{0} - v_{0})_+\leq 0,\end{equation}
	which implies that $u \le v$ a.e. in $\Omega\times (0, T)$. Since $\Omega$ was arbitrary, the comparison principle is proved. 
\end{proof}
	
\section{Estimates of the auxiliary problem}\label{secaux}
Let $\Omega$ be a precompact open subset of $M$.
Consider for $N>1$ the following auxiliary Cauchy problem:

\begin{equation}
	\left\{ 
	\begin{array}{ll}
		\partial _{t}u=\text{div}\left(A(u, \nabla u)\right) &\text{in}~\Omega\times (0, T), \\ 
		u=\frac{1}{N}, & \text{on}~\partial \Omega\times (0, T),\\u(x, 0)=u_{0}(x)+\frac{1}{N}& \text{in}~\Omega,
	\end{array}%
	\right.  \label{dtvaux}
\end{equation}%
where $A(u, \nabla u)=q^{p-1}\chi(u)^{(q-1)(p-1)}|\nabla u|^{p-2}\nabla u$ and $\chi(u)=\min\left(N, \max\left(u, \frac{1}{N}\right)\right)$.

The Cauchy problem (\ref{dtvaux}) is understood in the weak sense as follows: 

\begin{definition}\label{defapp}
We say that $u=u(x, t)$ is a \textit{weak subsolution} of (\ref{dtvaux}) in $\Omega\times (0, T)$, if 
\begin{equation}  \label{defvonsoluqaux}
	u\in C\left([0, T]; L^{1}(\Omega)\right) \text{\ and\ }
	u-\frac{1}{N}\in L^{p}\left((0, T); W_{0}^{1, p}(\Omega)\right)
\end{equation}
and for all \textit{test functions} 
\begin{equation}  \label{defvontestsoluqaux}
	\psi\in W^{1, \infty}\left((0, T);L^{\infty}(\Omega)\right)
	\cap L^{p}\left((0, T); W_{0}^{1, p}(\Omega)\right),
\end{equation}
and for all $t_{1}, t_{2}\in [0, T]$ with $t_{1}<t_{2}$, we have 
\begin{equation}  \label{defvonweaksolqaux}
	\left[\int_{\Omega}{u\psi}\right]_{t_{1}}^{t_{2}}+\int_{t_{1}}^{t_{2}}{%
		\int_{\Omega}{-u\partial_{t}\psi+\langle A(u, \nabla u),
			\nabla \psi\rangle}}\leq 0.
\end{equation}
Again, the notion of a \textit{weak supersolution} and \textit{weak solution} is defined analogously.
\end{definition} 

\begin{lemma}\label{lemexist}
Suppose that $u_{0}\in L^{1}(\Omega)\cap L^{\infty}(\Omega)$. Then, for any $N>1$, there exists a weak solution $u$ of the auxiliary problem (\ref{dtvaux}).
\end{lemma}

\begin{proof}
In this proof, we will use a modification of Galerkin's method given in \cite{ladyzhenskaya1968linear} (p. 466) (see also \cite{dekkers2005finite} and Lemma 5.2 in \cite{bogelein2021existence}).	
Consider the vector field $$\widetilde{A}(w, \xi)=A\left(\frac{1}{N}+w, \xi\right)=q^{p-1}\chi\left(\frac{1}{N}+w\right)^{(q-1)(p-1)}|\xi|^{p-2}\xi.$$ We will show that $w\in C\left([0, T]; L^{2}(\Omega)\right)$  and $w\in L^{p}\left((0, T); W_{0}^{1, p}(\Omega)\right)$
and for all $t_{1}, t_{2}\in [0, T]$ with $t_{1}<t_{2}$,
\begin{equation}\label{solutionw}
	\left[\int_{\Omega}{w\psi}\right]_{t_{1}}^{t_{2}}+\int_{t_{1}}^{t_{2}}{%
		\int_{\Omega}{-w\partial_{t}\psi+\langle \widetilde{A}(w, \nabla w),
			\nabla \psi\rangle}}= 0
\end{equation}
for all $\psi\in W^{1, 2}\left((0, T);L^{2}(\Omega)\right)
	\cap L^{p}\left((0, T); W_{0}^{1, p}(\Omega)\right)$
and $w(x, 0)=u_{0}(x)$. 
This will imply that $u=\frac{1}{N}+w$ is a weak solution of the auxiliary problem (\ref{dtvaux}). 

Denote $V=L^{2}(\Omega)\cap W_{0}^{1, p}(\Omega)$ and let $\{\phi_{m}\}_{m=1}^{\infty}$ be a complete basis of $V$ such that $(\phi_{m}, \phi_{n})_{L^{2}(\Omega)}=\delta_{m n}$. Let us show that, for any $n$, there exists coefficients $c_{m}^{n}\in C^{1}([0, T])$ such that the function $$w_{n}(x, t)=\sum_{m=1}^{n}c_{m}^{n}(t)\phi_{m}(x)$$ is a solution to the problem \begin{equation}\label{finitegalerkin}\int_{\Omega}\frac{d}{dt}w_{n}(t)\phi_{m}+\int_{\Omega}\langle \widetilde{A}(w_{n}, \nabla w_{n}), \nabla \phi_{m}\rangle=0, \qquad m=1, \ldots, n\end{equation} and $w_{n}(\cdot, 0)=u_{0}$ a.e. in $\Omega$. To prove this, it is sufficient to have $$\frac{d}{dt}c_{m}^{n}(t)=-\int_{\Omega}\langle \widetilde{A}(w_{n}, \nabla w_{n}), \nabla \phi_{m}\rangle, \qquad m=1, \ldots, n$$ and \begin{equation}\label{initialode}c_{m}^{n}(0)=\int_{\Omega}u_{0}\phi_{m}, \qquad m=1, \ldots, n\end{equation} since then, for $m=1, \ldots, n$, $$\int_{\Omega}\frac{d}{dt}w_{n}(t)\phi_{m}=\sum_{l=1}^{n}\frac{d}{dt}c_{l}^{n}(t)\int_{\Omega}\phi_{l}\phi_{m}=\frac{d}{dt}c_{m}^{n}=-\int_{\Omega}\langle \widetilde{A}(w_{n}, \nabla w_{n}), \nabla \phi_{m}\rangle$$ and $$w_{n}(x, 0)=\sum_{m=1}^{n}c_{m}^{n}(0)\phi_{m}(x)=\sum_{m=1}^{n}\int_{\Omega}u_{0}\phi_{m}\phi_{m}(x)=u_{0}(x).$$
We have $$\langle \widetilde{A}(w, \xi), \xi\rangle=q^{p-1}\chi\left(\frac{1}{N}+w\right)^{(q-1)(p-1)}|\xi|^{p}\geq c|\xi|^{p},$$ where $c$ depends on $p, q$ and $N$. Thus, multiplying (\ref{finitegalerkin}) by $c_{m}^{n}(t)$, summing over $m$ and integrating over $[0, t]$, where $t\in [0, T]$, yields \begin{equation}\label{uniformbounded}\frac{1}{2}||w_{n}(t)||_{L^{2}(\Omega)}^{2}+C\int_{0}^{t}\int_{\Omega}|\nabla w_{n}|^{p}\leq \frac{1}{2}||w_{n}(0)||_{L^{2}(\Omega)}^{2}= \frac{1}{2}||u_{0}||_{L^{2}(\Omega)}^{2}.
\end{equation} Therefore, $$\sum_{m=1}^{n}(c_{m}^{n}(t))^{2}=||w_{n}(t)||_{L^{2}(\Omega)}^{2}\leq ||u_{0}||_{L^{2}(\Omega)}^{2},$$ so that the bound is uniform in $n$ and $t$. Hence, we can solve the ODE in the whole interval $[0, T]$ and obtain coefficients $c_{m}^{n}$ satisfying the above conditions. 
Note that $$|\widetilde{A}(w, \xi)|=q^{p-1}\chi\left(\frac{1}{N}+w\right)^{(q-1)(p-1)}|\xi|^{p-1}\leq C|\xi|^{p-1},$$
where $C$ depends on $p, q$ and $N$.
Therefore, using also (\ref{uniformbounded}), we obtain \begin{align*}
|c_{m}^{n}(t)-c_{m}^{n}(t+h)|&=\left| \int_{\Omega}(w_{n}(\cdot, t)-w_{n}(\cdot, t+h))\phi_{m}\right|\\&\leq \int_{t}^{t+h}\int_{\Omega}|\langle \widetilde{A}(w_{n}, \nabla w_{n}), \nabla \phi_{m}\rangle|\\&\leq C\left(\int_{t}^{t+h}\int_{\Omega}|\nabla w_{n}|^{p}\right)^{\frac{p-1}{p}}h^{\frac{1}{p}}||\nabla \phi_{m}||_{L^{p}(\Omega)}^{p}\\&\leq C_{m}||u_{0}||_{L^{2}(\Omega)}^{2\frac{(p-1)}{p}}h^{\frac{1}{p}}.
\end{align*}
Hence, the coefficients $c_{m}^{n}(t)=\int_{\Omega}w_{n}(\cdot, t)\phi_{m}$ are uniformly bounded and equicontinuous. From the Arzelà-Ascoli theorem and a diagonal argument, we obtain, for any $m$, a uniform limit $c_{m}$ of $c_{m}^{n}$.

Let us show that $w(x, t)=\sum_{m=1}^{\infty}c_{m}(t)\phi_{m}(x)$ is a solution of the problem from above. We have $$||w(\cdot, t)||_{L^{2}(\Omega)}^{2}=\lim_{n\to \infty}\sum_{m=1}^{n}(c_{m}(t))^{2}\leq \lim_{n\to \infty} \sum_{m=1}^{n}|c_{m}(t)-c_{m}^{s}(t)|^{2}+||w_{s}(t)||_{L^{2}(\Omega)}^{2}\leq C||u_{0}||_{L^{2}(\Omega)}^{2},$$ for large enough $s$, since $c_{m}^{s}(t)\to c_{m}(t)$ uniformly in $t$. Therefore, $w(\cdot, t)\in L^{2}(\Omega)$ and in particular, $w\in C([0, T], L^{2}(\Omega))$. Similarly, for any $f\in L^{2}(\Omega)$,  \begin{align*}
(w_{n}(\cdot, t)-w(\cdot, t), f)&=\int_{\Omega}\sum_{m=1}^{\infty}\left(\int_{\Omega}f\phi_{m}\right)(w_{n}(\cdot, t)-w(\cdot, t))\phi_{m}\\&=\sum_{m=1}^{s}\left(\int_{\Omega}f\phi_{m}\right)(c_{m}^{n}(t)-c_{m}(t))+\sum_{m=s+1}^{\infty}\left(\int_{\Omega}f\phi_{m}\right)(c_{m}^{n}(t)-c_{m}(t)),
\end{align*} which tends to $0$ as $n\to \infty$. Thus, $w_{n}(\cdot, t)$ converges weakly to $w(\cdot, t)$ in $L^{2}(\Omega)$. 
Using the local Poincaré inequality applied to $w_{n}(\cdot, t)\in W_{0}^{1, p}(\Omega)$, \begin{equation}\label{localpoincineq}\int_{\Omega}w_{n}^{p}\leq C(\Omega)\int_{\Omega}|\nabla w_{n}|^{p},\end{equation} where $ C(\Omega)$ depends on the local geometry of $\Omega$, we deduce from (\ref{uniformbounded}) that \begin{equation}\label{localpoinc}\int_{0}^{T}||w_{n}(s)||_{W^{1, p}(\Omega)}^{p}ds\leq C(\Omega)||u_{0}||_{L^{2}(\Omega)}^{2} .\end{equation}
This implies that $w\in L^{p}\left((0, T); W_{0}^{1, p}(\Omega)\right)$.

From (\ref{initialode}) we also have $$w(x, 0)=\sum_{m=1}^{\infty}c_{m}(0)\phi_{m}(x)=\lim_{n\to \infty}\sum_{m=1}^{n}\phi_{m}(x)\int_{\Omega}u_{0}\phi_{m}=u_{0}(x)\quad \textnormal{for a.e.}~x~\textnormal{in}~ \Omega.$$
Hence, it remains to show that $w$ satisfies (\ref{solutionw}). From (\ref{finitegalerkin}), we have
\begin{equation}\label{approxsolupr}\left[\int_{\Omega}{w_{n}\psi}\right]_{t_{1}}^{t_{2}}+\int_{t_{1}}^{t_{2}}{%
\int_{\Omega}{-w_{n}\partial_{t}\psi+\langle \widetilde{A}(w_{n}, \nabla w_{n}),	\nabla \psi\rangle}}= 0,\end{equation} for any function $\psi(x, t)=\sum_{m=1}^{n}d_{m}(t)\phi_{m}(x)$, where the $d_{m}$ are continuous with weak derivatives in $L^{2}([0, T])$. The set of the union over all $n$ of these functions is dense in the space $W^{1, 2}\left((0, T);L^{2}(\Omega)\right)
\cap L^{p}\left((0, T); W_{0}^{1, p}(\Omega)\right)$, whence it is sufficient to show (\ref{solutionw}) for all such $\psi$. From the weak $L^{2}$-convergence of $w_{n}$ to $w$, we deduce $$\lim_{n\to \infty}\left[\int_{\Omega}{w_{n}\psi}\right]_{t_{1}}^{t_{2}}\to \left[\int_{\Omega}{w\psi}\right]_{t_{1}}^{t_{2}}$$ and $$\lim_{n\to \infty}\int_{t_{1}}^{t_{2}}{\int_{\Omega}{w_{n}\partial_{t}\psi}}=\int_{t_{1}}^{t_{2}}{\int_{\Omega}{w\partial_{t}\psi}}.$$
Therefore, we are left to prove that \begin{equation}\label{wsatAlim}\lim_{n\to \infty}\int_{t_{1}}^{t_{2}}{\int_{\Omega}{\langle \widetilde{A}(w_{n}, \nabla w_{n}),	\nabla \psi\rangle}}=\int_{t_{1}}^{t_{2}}{\int_{\Omega}{\langle \widetilde{A}(w, \nabla w),	\nabla \psi\rangle}}.\end{equation}

Let us first show that $\nabla w_{n}\to \nabla w$ in $L^{p}(\Omega\times [t_{1}, t_{2}])$.
It follows from Lemma \ref{plapmon} that $$\langle \widetilde{A}(w_{n}, \nabla w_{n})-\widetilde{A}(w, \nabla w), \nabla w_{n}-\nabla w\rangle\geq c\left\{ 
\begin{array}{ll}
	|\nabla w_{n}-\nabla w|^{p} &\text{if}~p\geq2, \\ 
	W_{n}^{p-2}|\nabla w_{n}-\nabla w|^{2}, & \text{if}~p<2,
\end{array}%
\right.$$ where $W_{n}=\left(|\nabla w_{n}|^{2}+|\nabla w|^{2}\right)^{\frac{1}{2}}$ and $c$ depends on $p, q$ and $N$. Let us denote $Q=B\times[t_{1}, t_{2}]$. In the case $p<2$ we obtain by Hölder's inequality, \begin{align*}
	\int_{Q}|\nabla w_{n}-\nabla w|^{p}&\leq \int_{Q\cap \{\nabla w_{n}\ne\nabla w\}}|\nabla w_{n}-\nabla w|^{p}W_{n}^{\frac{p(p-2)}{2}}W_{n}^{\frac{p(2-p)}{2}}\\&\leq \left(\int_{Q\cap \{\nabla w_{n}\ne\nabla w\}}|\nabla w_{n}-\nabla w|^{2}W_{n}^{p-2}\right)^{p/2}\left(\int_{Q}W_{n}^{p}\right)^{(2-p)/2}\\&\leq C\left(\int_{Q\cap \{\nabla w_{n}\ne\nabla w\}}\langle \widetilde{A}(w_{n}, \nabla w_{n})- \widetilde{A}(w, \nabla w), \nabla w_{n}-\nabla w\rangle\right)^{p/2}\\&\quad\times\left(\int_{Q}|\nabla w_{n}|^{p}+|\nabla w|^{p}\right)^{(2-p)/2}.
\end{align*}
By (\ref{localpoinc}), the latter integral is bounded by some constant independant of $n$. Therefore, \begin{align}
	\frac{1}{C}\left(\int_{Q}|\nabla w_{n}-\nabla w|^{p}\right)^{\gamma}&\nonumber\leq \int_{Q}\langle \widetilde{A}(w_{n}, \nabla w_{n})- \widetilde{A}(w, \nabla w), \nabla w_{n}-\nabla w\rangle.
\end{align} where $\gamma=\min(1, 2/p)$ and $C$ does not depend on $n$. Hence, we have in all cases \begin{align*}
\left(\int_{Q}|\nabla w_{n}-\nabla w|^{p}\right)^{\gamma}&\leq C\int_{Q}\langle \widetilde{A}(w_{n}, \nabla w_{n})-\widetilde{A}(w, \nabla w), \nabla w_{n}-\nabla w\rangle\\&=C\int_{Q}\langle \widetilde{A}(w_{n}, \nabla w_{n}),\nabla w_{n}\rangle-\langle \widetilde{A}(w_{n}, \nabla w_{n}),\nabla w\rangle-\langle \widetilde{A}(w, \nabla w), \nabla w_{n}-\nabla w\rangle\\&\leq C\left(\frac{1}{2}||w_{n}(\cdot, 0)||_{L^{2}(\Omega)}^{2}-\frac{1}{2}||w_{n}(\cdot, t_{2})||_{L^{2}(\Omega)}^{2}\right)\\&\quad-C\int_{Q}\langle \widetilde{A}(w_{n}, \nabla w_{n}),\nabla w\rangle+\langle \widetilde{A}(w, \nabla w), \nabla w_{n}-\nabla w\rangle,
\end{align*}
where we used a similar argument as in (\ref{uniformbounded}).
Note that (\ref{uniformbounded}) also implies that $\nabla w_{n}\to \nabla w$ weakly in $L^{p}(Q)$ and that  $\widetilde{A}(w_{n}, \nabla w_{n})$ converges weakly in $L^{\frac{p}{p-1}}(Q)$ to some $\xi:Q\to TM$, which yields \begin{equation}\label{gradlpconvw}\lim_{n\to \infty}\left(\int_{Q}|\nabla w_{n}-\nabla w|^{p}\right)^{\gamma}\leq C\lim_{n\to \infty}\left(\frac{1}{2}||w_{n}(\cdot, 0)||_{L^{2}(\Omega)}^{2}-\frac{1}{2}||w_{n}(\cdot, t_{2})||_{L^{2}(\Omega)}^{2}\right)-\int_{Q}\langle \xi, \nabla w\rangle.
\end{equation} It follows from (\ref{approxsolupr}) that, for all $0<h<T$ and for all $\tau\in [0, T]$, \begin{equation}\label{steklovxi}\int_{0}^{\tau}\int_{\Omega}{(\partial_{t} w^{h})\psi+\langle\xi^{h}, \nabla \psi\rangle }= 0,\end{equation} where $w^{h}$ is defined by (\ref{}). Testing with $\psi=w^{h}$ in (\ref{steklovxi}), we obtain $$\int_{Q}\langle \xi^{h}, \nabla w\rangle= \frac{1}{2}||w^{h}(\cdot, 0)||_{L^{2}(\Omega)}^{2}-\frac{1}{2}||w^{h}(\cdot, t_{2})||_{L^{2}(\Omega)}^{2} ,$$ whence sending $h\to 0$ yields (c.f. Lemma 2.2 in \cite{dekkers2005finite}) $$\int_{Q}\langle \xi, \nabla w\rangle= \frac{1}{2}||u_{0}||_{L^{2}(\Omega)}^{2}-\frac{1}{2}||w(\cdot, t_{2})||_{L^{2}(\Omega)}^{2}.$$ Let us also notice that $w_{n}(\cdot, 0)\to u_{0}$ in $L^{2}(\Omega)$ as $$||w_{n}(\cdot, 0)-u_{0}||_{L^{2}(\Omega)}^{2}=\sum_{m=n+1}^{\infty}|c_{m}(0)|^{2}\to 0 \quad \textnormal{as}~n\to \infty,$$ and by the weak convergence of $w_{n}(\cdot, t)\to w(\cdot, t)$ in $L^{2}(\Omega)$, $$\liminf_{n\to \infty}||w_{n}(\cdot, t_{2})||_{L^{2}(\Omega)}\geq ||w(\cdot, t_{2})||_{L^{2}(\Omega)}^{2}.$$
Combining these estimates with (\ref{gradlpconvw}), we deduce $$\liminf_{n\to \infty}\left(\int_{Q}|\nabla w_{n}-\nabla w|^{p}\right)^{\gamma}=0,$$ so that $\nabla w_{n}\to \nabla w$ in $L^{p}(Q)$.

Let us apply the following inequality, which holds for all $\xi, \eta\in T_{x}M$, \begin{equation}\label{inequalitygrad}\left||\xi|^{p-2}\xi-|\eta|^{p-2}\eta\right|\leq C(|\eta|+|\xi-\eta|)^{p-2}|\xi-\eta|\end{equation} (see Lemma 3.6 in \cite{bogelein2021existence}).
Hence, we have \begin{align*}|\langle \widetilde{A}(w_{n}, \nabla w_{n})-\widetilde{A}(w, \nabla w),	\nabla \psi\rangle|&\leq C\left||\nabla w_{n}|^{p-2}\nabla w_{n}-|\nabla w|^{p-2}\nabla w\right||\nabla \psi|\\&\leq C(|\nabla w|+|\nabla w_{n}-\nabla w|)^{p-2}|\nabla w_{n}-\nabla w||\nabla \psi|\\&\leq (C_{1}|\nabla w|^{p-2}|\nabla w_{n}-\nabla w|+C|\nabla w_{n}-\nabla w|^{p-1})|\nabla \psi|,\end{align*} where $C$ depends on $p, q$ and $N$ and $C_{1}=0$ if $p<2$. In the case $p<2$, it follows that $$\int_{Q}|\langle \widetilde{A}(w_{n}, \nabla w_{n})-\widetilde{A}(w, \nabla w),	\nabla \psi\rangle|\leq C\int_{Q}|\nabla w_{n}-\nabla w|^{p-1}|\nabla \psi|\leq C||\nabla w_{n}-\nabla w||_{L^{p}}||\nabla \psi||_{L^{p}}.$$ On the other hand, if $p\geq 2$, we also use \begin{align*}C_{1}\int_{Q}|\nabla w|^{p-2}|\nabla w_{n}-\nabla w||\nabla \psi|&\leq C_{1}\left(\int_{Q}|\nabla \psi|^{\frac{p}{p-1}}|\nabla w|^{\frac{p(p-2)}{p-1}}\right)^{\frac{p-1}{p}}||\nabla w_{n}-\nabla w||_{L^{p}}\\&\leq C_{1}||\nabla \psi||_{L^{p}}||\nabla w||_{L^{p}}^{p-2}||\nabla w_{n}-\nabla w||_{L^{p}}.
\end{align*} Therefore, we conclude that $$\lim_{n\to \infty}\int_{Q}{\langle \widetilde{A}(w_{n}, \nabla w_{n}),	\nabla \psi\rangle}=\int_{Q}{\langle \widetilde{A}(w, \nabla w),	\nabla \psi\rangle},$$  which proves (\ref{wsatAlim}) and finishes the proof. 
\end{proof}

\begin{corollary}\label{loweru}
	Let $u$ be a weak solution of (\ref{dtvaux}) and assume that $u_{0}\geq 0$. Then $u\geq \frac{1}{N}$ a.e. in $\Omega\times (0, T)$.
\end{corollary}

\begin{proof}
	Applying (\ref{forcomparison}) with subsolution $v=\frac{1}{N}$ and $v_{0}=0$, we get $u\geq v=\frac{1}{N}$ a.e. in $\Omega\times (0, T)$.
\end{proof}

\begin{lemma}\label{lipschitzapp}
Let $u$ be a weak solution of (\ref{dtvaux}) in $\Omega\times (0, T)$ and $g$ be a Lipschitz function in $\mathbb{R}_{+}$ such that $g(1/N)=0$. Then, for all $t\in [0, T]$, \begin{equation}\label{withlipschitz}
		\left[\int_{\Omega}{G(u)}\right]_{0}^{t}+\int_{0}^{t}{%
			\int_{\Omega}{\langle A(u, \nabla u),\nabla u\rangle g^{\prime}(u)}}= 0,
\end{equation}
where $G(u)=\int_{0}^{u}{g(s)ds}$.
\end{lemma}

\begin{proof}
We have, for all $0<h<T$ and for all $\tau\in [0, T]$, $$\int_{0}^{\tau}\int_{\Omega}{(\partial_{t} u^{h})\psi+\langle[A(u, \nabla u)]^{h}, \nabla \psi\rangle }= 0,$$ where $u^{h}$ is defined as in (\ref{defclassicstekl}). Choosing $\psi=g(u^{h})$, which is admissible as a testfunction since $g(1/N)=0$, we obtain $$\left[\int_{\Omega}{G(u^{h})}\right]_{0}^{\tau}+\int_{0}^{\tau}\int_{\Omega}{\langle[A(u, \nabla u)]^{h}, \nabla u^{h}\rangle g^{\prime}(u^{h}) }= 0.$$ Hence, letting $h\to 0$, we get (\ref{withlipschitz}).
\end{proof}

\begin{lemma}\label{upperu}
Let $u$ be a weak solution of (\ref{dtvaux}) in $\Omega\times (0, T)$. Then, for all $N\geq N_{0}(u_{0})$, $$||u||_{L^{\infty}(\Omega\times (0, T))}\leq C,$$ where $C$ depends on $p, q$ and $u_{0}$, but not on $N$.
\end{lemma}

\begin{proof}
Let us apply Lemma \ref{lipschitzapp} with $g(u)=\left(\chi(u)-||u_{0}||_{L^{\infty}(\Omega)}-1/N\right)_{+}$. Then it follows from (\ref{withlipschitz}) with $G(u)\geq \frac{1}{2}g(u)^{2}$ that, for all $t\in [0, T]$,
$$\left[\int_{\Omega}{g(u)^{2}}\right]_{0}^{t}+2\int_{0}^{t}{\int_{\Omega}{\langle A(u, \nabla u),\nabla u\rangle g^{\prime}(u)}}\leq 0.$$ Noticing that $$\int_{0}^{t}{\int_{\Omega}{\langle A(u, \nabla u),\nabla u\rangle g^{\prime}(u)}}\geq c\int_{0}^{t}{\int_{\Omega}|\nabla u|^{p}g^{\prime}(u)}\geq0$$ and $\int_{\Omega}{g(u)(\cdot, 0)^{2}}=0$ we obtain $$\sup_{t\in (0, T)}\int_{\Omega}{g(u)^{2}}\leq 0.$$ 
Therefore, $$\chi(u)\leq ||u_{0}||_{L^{\infty}(\Omega)}+1/N.$$ Now choosing $N>N_{0}(u_{0})$ large enough, we deduce $u\leq C$, with a constant $C$ depending on $p, q$ and $u_{0}$ but not on $N$.
\end{proof}

\subsection{Limiting procedure for solutions of the auxiliary problem}

Since $\Omega$ is precompact, it holds that $W_{0}^{1, p}(\Omega)$ is compactly embedded in $L^{p}(\Omega)$. Hence, we have the following result.
\begin{lemma}[Corollary 4.6(i) in \cite{bogelein2021existence}]\label{lemmaforstrong}
	Let $f_{n}$ be a bounded sequence in $L^{\theta}((0, T), L^{\frac{p}{p-1}}(\Omega))$, where $\theta=\max(1, pq)$, so that $(f_{n})^{q}$ is bounded in $L^{p}((0, T), W^{1, p}(\Omega))$. If \begin{equation}\label{condforstrong}||f_{n}(t+h, x)-f_{n}(x)||_{L^{\theta}\left((0, T-h), \left(W_{0}^{1, p}(\Omega)\right)^{\prime}\right)}\to 0,~\textnormal{uniformly as}~h\to 0,\end{equation} then $(f_{n})^{q}$ is relatively compact in $L^{p}(\Omega\times(0, T))$.		
\end{lemma}

\begin{lemma}\label{weakstrongconv}
Let $u_{N}$ be the solution to the Cauchy problem (\ref{dtvaux}) in $\Omega\times (0, T)$ from Lemma \ref{lemexist}. Then there exists a subsequence $\{N_{j}\}_{j}$ and a function $u\in L^{\infty}(\Omega\times (0, T))$ with $u^{q}\in L^{p}((0, T), W_{0}^{1, p}(\Omega))$ such that \begin{equation}\label{weakcon}u_{N_{j}}^{q}\to u^{q}\quad\textnormal{weakly in}~ L^{p}((0, T), W^{1, p}(\Omega))\end{equation} and \begin{equation}\label{strongcon}u_{N_{j}}\to u\quad\textnormal{strongly in}~ L^{r}(\Omega\times (0, T))~\textnormal{for any}~r\geq 1~ \textnormal{and a.e. in}~\Omega\times (0, T).\end{equation}
\end{lemma}

\begin{proof}
Combining Corollary \ref{loweru} and Lemma \ref{upperu}, we get for $N>\max\{N_{0}, C\}$ that $\frac{1}{N}<u_{N}<C<N$ and thus \begin{equation}\label{largeNLeib}A(u_{N}, \nabla u_{N})=q^{p-1}|\nabla u_{N}^{q}|^{p-2}\nabla u_{N}^{q}.\end{equation} 
Applying Lemma \ref{lipschitzapp} with $g(s)=\left(s^{q}-N^{-q}\right)_{+}$, we deduce 
with $G(s)\geq\frac{1}{q+1}g(s)^{\frac{q+1}{q}}$, since $u_{N}>\frac{1}{N}$ and using (\ref{largeNLeib}), \begin{equation*}\frac{1}{q+1}\left[\int_{\Omega}{g(u_{N})^{\frac{q+1}{q}}}\right]_{0}^{T}+\int_{Q}{|\nabla u_{N}^{q}|^{p}}\leq \left[\int_{\Omega}{G(u_{N})}\right]_{0}^{T}+\int_{0}^{T}{\int_{\Omega}{\langle A(u_{N}, \nabla u_{N}),\nabla u_{N}\rangle g^{\prime}(u_{N})}}= 0,\end{equation*} where $Q=\Omega\times (0, T)$. We have for $N>N(u_{0})$,  $$\int_{\Omega}{g(u_{N}(\cdot, 0))^{\frac{q+1}{q}}}\leq C\int_{\Omega}u_{0}^{q+1},$$ and therefore,  \begin{equation}\label{uppergrad}\int_{Q}|\nabla u_{N}^{q}|^{p}\leq C\int_{\Omega}u_{0}^{q+1}.\end{equation}
On the other hand, using the local Poincaré inequality (\ref{localpoincineq}) applied to $(u_{N}^{q}-N^{-q})_{+}(\cdot, t)\in W_{0}^{1, p}(\Omega)$, we deduce from (\ref{uppergrad}) that \begin{equation*}\int_{Q}(u_{N}^{q}-N^{-q})_{+}^{p}\leq C(\Omega)\int_{\Omega}u_{0}^{q+1},\end{equation*} where $C(\Omega)$ depends on the local geometry of $\Omega$.
Therefore, the sequence $(u_{N}^{q}-N^{-q})$ is a bounded sequence in $L^{p}((0, T), W_{0}^{1, p}(\Omega))$. Hence, there is a subsequence $\{N_{j}\}_{j}$ with $N_{j}\to \infty$ and a function $w\in L^{p}((0, T), W_{0}^{1, p}(\Omega))$ such that $(u_{N_{j}}^{q}-N_{j}^{-q})\to w$ weakly in $L^{p}((0, T), W^{1, p}(\Omega))$. This yields that \begin{equation}\label{weakconsub}u_{N_{j}}^{q}\to w\quad\textnormal{weakly in}~ L^{p}((0, T), W^{1, p}(\Omega)).\end{equation}

Let us now show the strong convergence of the sequence $u_{N_{j}}$ in $L^{p}((0, T), L^{p}(\Omega))=L^{p}(\Omega\times (0, T))$. For that, we want to apply Lemma \ref{lemmaforstrong}. For that purpose, let us apply equation (\ref{defvonweaksolqaux}) and (\ref{largeNLeib}) with test function $\psi(x, t)=\psi(x)\in W_{0}^{1, p}(\Omega)$. Hence, we have $$ \left|\int_{\Omega}\left(u_{N}(t+h)-u_{N}(t)\right)\psi\right|\leq \int_{t}^{t+h} \int_{\Omega}|\nabla u_{N}^{q}|^{p-1}|\nabla \psi|.$$ By Hölder's inequality and (\ref{uppergrad}), this implies that $$\left|\int_{\Omega}\left(u_{N}(t+h)-u_{N}(t)\right)\psi\right|\leq Ch^{1/p}||\psi||_{W^{1, p}(\Omega)}.$$ Denoting with $\langle \cdot, \cdot\rangle_{\left(W_{0}^{1, p}(\Omega)\right)^{\prime}\times W_{0}^{1, p}(\Omega)}$ the dual pairing of $\left(W_{0}^{1, p}(\Omega)\right)^{\prime}$ and $W_{0}^{1, p}(\Omega)$ and regarding $u_{N}(t)$ as an element of the dual space $\left(W_{0}^{1, p}(\Omega)\right)^{\prime}$ for $t\in [0, T]$ we obtain \begin{align*}||u_{N}(t+h)-u_{N}(t)||_{\left(W_{0}^{1, p}(\Omega)\right)^{\prime}}&=\sup_{0\ne\psi\in W_{0}^{1, p}(\Omega)}\frac{\left|\langle u_{N}(t+h)-u_{N}(t), \psi\rangle_{\left(W_{0}^{1, p}(\Omega)\right)^{\prime}\times W_{0}^{1, p}(\Omega)} \right|}{||\psi||_{W^{1, p}(\Omega)}}\\&\leq\sup_{0\ne\psi\in W_{0}^{1, p}(\Omega)}\frac{Ch^{1/p}||\psi||_{W^{1, p}(\Omega)}}{||\psi||_{W^{1, p}(\Omega)}}=Ch^{1/p}.
\end{align*} This proves (\ref{condforstrong}) for $u_{N}$. Further, we know from Lemma \ref{upperu} that $u_{N}$ is a bounded sequence in $L^{\infty}(Q)$ so that it is bounded in $L^{\theta}((0, T), L^{\frac{p}{p-1}}(\Omega))$. From (\ref{uppergrad}) we know that $u_{N}^{q}$ is a bounded sequence in $L^{p}((0, T), W^{1, p}(\Omega))$. This yields the relative compactness of the sequence $\{u_{N_{j}}^{q}\}_{j}$ in $L^{p}(Q)$. Together with (\ref{weakconsub}) this implies that \begin{equation}\label{strongconsub}u_{N_{j}}^{q}\to w\quad\textnormal{strongly in}~ L^{p}(Q).\end{equation}
By passing to yet another subsequence, also denoted by $u_{N_{j}}$, we also obtain that $u_{N_{j}}^{q}\to w$ a.e. in $Q$ and that $w$ is bounded in $Q$ since the sequence $u_{N_{j}}$ is uniformly bounded. If we now define $u:=w^{1/q}$, we get from (\ref{weakconsub}) and (\ref{strongconsub}) that $u_{N_{j}}^{q}\to u^{q}$ weakly in $L^{p}((0, T), W^{1, p}(\Omega))$ and that $u_{N_{j}}^{q}\to u^{q}$ in $L^{p}(Q)$. Using the latter convergence of the sequence $u_{N_{j}}^{q}$ and its uniform boundedness, we also obtain the convergence of $u_{N_{j}}^{q}$ to $u^{q}$ in $L^{r}(Q)$ for any $r\geq 1$. This finishes the proof.
\end{proof}

\begin{lemma}\label{convofgrad}
Under the assumptions of Lemma \ref{weakstrongconv} we have for any subinterval $I$ compactly contained in $(0, T)$, \begin{equation}\label{convgrad}\nabla u_{N_{j}}^{q}\to \nabla u^{q}\quad \textnormal{strongly in}~L^{p}(\Omega\times I).
\end{equation}
\end{lemma}

\begin{proof}
It follows from Lemma \ref{plapmon} that $$\langle |\nabla u_{N}^{q}|^{p-2}\nabla u_{N}^{q}- |\nabla u^{q}|^{p-2}\nabla u^{q}, \nabla u_{N}^{q}-\nabla u^{q}\rangle\geq \left\{ 
\begin{array}{ll}
	c|\nabla u_{N}^{q}-\nabla u^{q}|^{p} &\text{if}~p\geq2, \\ 
	cV_{N}^{p-2}|\nabla u_{N}^{q}-\nabla u^{q}|^{2}, & \text{if}~p<2,
\end{array}%
\right.$$ where $V_{N}=\left(|\nabla u_{N}^{q}|^{2}+|\nabla u^{q}|^{2}\right)^{\frac{1}{2}}$. Let $\zeta$ be a non-negative bounded Lipschitz function in $[0, T]$ such that $\zeta(0)=\zeta(T)=0$. Then we obtain in the case $p<2$ by Hölder's inequality, \begin{align*}
\int_{Q}\zeta|\nabla u_{N}^{q}-\nabla u^{q}|^{p}&\leq \int_{Q\cap \{\nabla u_{N}^{q}\ne\nabla u^{q}\}}\zeta|\nabla u_{N}^{q}-\nabla u^{q}|^{p}V_{N}^{\frac{p(p-2)}{2}}V_{N}^{\frac{p(2-p)}{2}}\\&\leq \left(\int_{Q\cap \{\nabla u_{N}^{q}\ne\nabla u^{q}\}}\zeta|\nabla u_{N}^{q}-\nabla u^{q}|^{2}V_{N}^{p-2}\right)^{p/2}\left(\int_{Q}\zeta V_{N}^{p}\right)^{(2-p)/2}\\&\leq C\left(\int_{Q\cap \{\nabla u_{N}^{q}\ne\nabla u^{q}\}}\zeta\langle |\nabla u_{N}^{q}|^{p-2}\nabla u_{N}^{q}- |\nabla u^{q}|^{p-2}\nabla u^{q}, \nabla u_{N}^{q}-\nabla u^{q}\rangle\right)^{p/2}\\&\quad\times\left(\int_{Q}|\nabla u_{N}^{q}|^{p}+|\nabla u^{q}|^{p}\right)^{(2-p)/2}.
\end{align*}
By (\ref{uppergrad}), the latter integral is bounded by some constant independant of $N$. Therefore, in all cases \begin{align}
\frac{1}{C}\left(\int_{Q}\zeta|\nabla u_{N}^{q}-\nabla u^{q}|^{p}\right)^{\gamma}&\nonumber\leq \int_{Q}\zeta\langle |\nabla u_{N}^{q}|^{p-2}\nabla u_{N}^{q}, |\nabla u^{q}|^{p-2}\nabla u^{q}- \nabla u_{N}^{q}-\nabla u^{q}\rangle\\&\nonumber=\int_{Q}\zeta\langle |\nabla u_{N}^{q}|^{p-2}\nabla u_{N}^{q}, \nabla u_{N}^{q}-\nabla u^{q}\rangle-\zeta\langle |\nabla u^{q}|^{p-2}\nabla u^{q}, \nabla u_{N}^{q}-\nabla u^{q}\rangle\\&\label{INIIN}=:I_{N}-II_{N},
\end{align} where $\gamma=\min(1, 2/p)$ and $C$ does not depend on $N$. 

Since $\zeta|\nabla u^{q}|^{p-2}\nabla u^{q}\in L^{\frac{p}{p-1}}(Q)$ and $u_{N_{j}}\to u$ weakly in $L^{p}((0, T), W_{0}^{1, p}(\Omega))$ for some subsequence $\{N_{j}\}_{j}$ by Lemma \ref{weakstrongconv}, we get $$\lim_{j\to \infty}II_{N_{j}}=\lim_{j\to \infty}\int_{Q}\zeta\langle |\nabla u^{q}|^{p-2}\nabla u^{q}, \nabla u_{N_{j}}^{q}-\nabla u^{q}\rangle=0.$$ 

Let us now write $$I_{N}=\int_{Q}\zeta\langle |\nabla u_{N}^{q}|^{p-2}\nabla u_{N}^{q}, \nabla u_{N}^{q}-\nabla (u^{q})_{h}\rangle+\int_{Q}\zeta\langle |\nabla u_{N}^{q}|^{p-2}\nabla u_{N}^{q}, \nabla (u^{q})_{h}-\nabla u^{q}\rangle=:I_{N}^{(1)}+I_{N}^{(2)},$$ where $(u^{q})_{h}$ is defined as in (\ref{defsteklov}). By Hölder's inequality and (\ref{uppergrad}), we have \begin{equation}\label{I2N}I_{N}^{(2)}\leq |||\nabla u_{N}^{q}|^{p-1}||_{L^{\frac{p}{p-1}}(Q)}||\nabla (u^{q})_{h}-\nabla u^{q}||_{L^{p}(Q)}\leq C||\nabla (u^{q})_{h}-\nabla u^{q}||_{L^{p}(Q)},\end{equation} with $C$ not depending on $N$.

Let $w\in L^{q+1}(Q)$ such that $w^{q}-N^{-q}\in L^{p}\left((0, T); W_{0}^{1,p}(\Omega)\right)$ and $\partial_{t}w^{q}\in L^{\frac{q+1}{q}}(Q)$.
Testing with $\psi=\zeta(w^{q}-(u_{N}^{q})_{h})$, where $w^{q}=w_{h, N}^{q}=(u^{q})_{h}+N^{-q}$, in (\ref{defvonweaksolq}), which is admissible since $u_{N}-\frac{1}{N}\in L^{p}\left((0, T); W_{0}^{1,p}(\Omega)\right)$, we obtain
similarly as in (\ref{laterusagefromcont}) in the proof of Lemma \ref{continuitylemma} that
\begin{align*}\int_{Q}\zeta^{\prime}b(u_{N}, w)-\left[\int_{\Omega} \zeta b(u_{N}, w)\right]_{0}^{T}=\int_{Q}\zeta(\partial_{t}w^{q}(u_{N}-w)+|\nabla u_{N}^{q}|^{p-2}\langle\nabla u_{N}^{q},\nabla u_{N}^{q}-\nabla w^{q}\rangle).\end{align*}
We get
$$I_{N}^{(1)}=-\int_{Q}\zeta\partial_{t}(u^{q})_{h}(u_{N}-w)+\int_{Q}\zeta^{\prime}b(u_{N}, w)-\left[\int_{\Omega}\zeta b(u_{N}, w)\right]_{0}^{T}.$$
Hence, using $\zeta(0)=\zeta(T)=0$ and (\ref{propmol}), \begin{align*}\limsup_{j\to\infty}I_{N_{j}}^{(1)}&=-\int_{Q}\zeta\partial_{t}(u^{q})_{h}(u-(u^{q})_{h}^{1/q})+\int_{Q}\zeta^{\prime}b(u, (u^{q})_{h}^{1/q})\\&=-\frac{1}{h}\int_{Q}(u^{q}-(u^{q})_{h})(u-(u^{q})_{h}^{1/q})+\int_{Q}\zeta^{\prime}b(u, (u^{q})_{h}^{1/q})\\&\leq\int_{Q}\zeta^{\prime}b(u, (u^{q})_{h}^{1/q}).\end{align*}
Combining this with (\ref{INIIN}) and (\ref{I2N}), we deduce \begin{align*}\limsup_{j\to\infty}\left(\int_{Q}\zeta|\nabla u_{N_{j}}^{q}-\nabla u^{q}|^{p}\right)^{\gamma}\leq C||\nabla (u^{q})_{h}-\nabla u^{q}||_{L^{p}(Q)}+\int_{Q}\zeta^{\prime}b(u, (u^{q})_{h}^{1/q}).\end{align*} Finally, by sending $h\to 0$, we conclude (\ref{convgrad}) from Lemma \ref{propstek}.
\end{proof}

\begin{lemma}\label{continuitylemmaap}
Let $u$ be the function obtained in Lemma \ref{weakstrongconv}. Then $u\in C([0, T], L^{q+1}(\Omega))$.
\end{lemma}

\begin{proof}
Let $u_{N_{j}}$ be the sequence of solutions of (\ref{dtvaux}) obtained in Lemma \ref{weakstrongconv}, which satisfies, for all $t_{1}, t_{2}\in [0, T]$ with $t_{1}<t_{2}$, $$\left[\int_{\Omega}{u_{N_{j}}\psi}\right]_{t_{1}}^{t_{2}}+\int_{t_{1}}^{t_{2}}{\int_{\Omega}{-u_{N_{j}}\partial_{t}\psi+|\nabla u_{N_{j}}^{q}|^{p-2}\langle\nabla u_{N_{j}}^{q},\nabla \psi\rangle}}= 0.$$
Choosing $\psi$ such that $\psi(0)=\psi(T)=0$ and $t_{1}=0$, $t_{2}=T$, we get by Lemma \ref{weakstrongconv}, $$\lim_{j\to \infty}\int_{0}^{T}{\int_{\Omega}{u_{N_{j}}\partial_{t}\psi}}=\int_{0}^{T}\int_{\Omega}{u\partial_{t}\psi}.$$ 
For the elliptic term $\mathcal{A}_{N_{j}}=|\nabla u_{N_{j}}^{q}|^{p-2}\nabla u_{N_{j}}^{q}-|\nabla u^{q}|^{p-2}\nabla u^{q}$, let us split the integration  $$\int_{0}^{T}\int_{\Omega}\langle\mathcal{A}_{N_{j}}, \nabla \psi\rangle=\int_{0}^{\tau}\int_{\Omega}\langle\mathcal{A}_{N_{j}}, \nabla \psi\rangle+\int_{\tau_{1}}^{\tau_{2}}\int_{\Omega}\langle\mathcal{A}_{N_{j}}, \nabla \psi\rangle+\int_{\tau_{2}}^{T}\int_{\Omega}\langle\mathcal{A}_{N_{j}}, \nabla \psi\rangle,$$ where $0<\tau_{1}<\tau_{2}<T$. Choosing $\tau_{1}$ small enough and $\tau_{2}$ large enough, we get from (\ref{uppergrad}) and Hölder's inequality, for any $\varepsilon>0$, $$\int_{0}^{\tau_{1}}\int_{\Omega}\left|\langle\mathcal{A}_{N_{j}}, \nabla \psi\rangle\right|<\varepsilon\quad \textnormal{and}\quad\int_{\tau_{2}}^{T}\int_{\Omega}\left|\langle\mathcal{A}_{N_{j}}, \nabla \psi\rangle\right|<\varepsilon.$$ Since $\nabla u_{N_{j}}^{q}\to\nabla u^{q}$ in $L^{p}(\Omega\times (\tau_{1}, \tau_{2}))$ by Lemma \ref{convofgrad} as $j\to \infty$, we can argue as in the proof of Lemma \ref{lemexist} and deduce that $$\lim_{j\to \infty}\int_{\tau_{1}}^{\tau_{2}}\int_{\Omega}\left|\langle\mathcal{A}_{N_{j}}, \nabla \psi\rangle\right|=0,$$ so that $$\int_{0}^{T}\int_{\Omega}{-u\partial_{t}\psi+|\nabla u^{q}|^{p-2}\langle\nabla u^{q},\nabla \psi\rangle}= 0.$$ Then the claim follows from Lemma \ref{continuitylemma}. 
\end{proof}

\begin{lemma}\label{lemmaexball}
Assume that $u_{0}\in L^{1}(\Omega)\cap L^{\infty}(\Omega)$ is non-negative. Then there exists a non-negative bounded weak solution $u$ of the Cauchy problem \begin{equation}
	\left\{ 
	\begin{array}{ll}
		\partial _{t}u=\Delta _{p}u^{q} &\text{in}~\Omega\times (0, T), \\ 
		u=0, & \text{on}~\partial \Omega\times (0, T),\\u(x, 0)=u_{0}(x)& \text{in}~\Omega,
	\end{array}%
	\right.  \label{Balldtv}
\end{equation} in the sense of Definition \ref{defweaksolu}.
\end{lemma}

\begin{proof}
Let $u$ be the function obtained in Lemma \ref{weakstrongconv}.
Let us first show that $u$ satisfies, for all $t_{1}, t_{2}\in [0, T]$ with $t_{1}<t_{2}$, \begin{equation}\label{usatsifies}
	\left[\int_{\Omega}{u\psi}\right]_{t_{1}}^{t_{2}}+\int_{t_{1}}^{t_{2}}{%
		\int_{\Omega}{-u\partial_{t}\psi+|\nabla u^{q}|^{p-2}\langle\nabla u^{q},
			\nabla \psi\rangle}}= 0,
\end{equation}
for all $\psi\in W^{1, \infty}\left((0, T);L^{\infty}(\Omega)\right)
\cap L^{p}\left((0, T); W_{0}^{1, p}(\Omega)\right)$.
For that, let $u_{N_{j}}$ be the sequence of solutions of (\ref{dtvaux}) obtained in Lemma \ref{weakstrongconv}, which satisfies $$\left[\int_{\Omega}{u_{N_{j}}\psi}\right]_{t_{1}}^{t_{2}}+\int_{t_{1}}^{t_{2}}{%
	\int_{\Omega}{-u_{N_{j}}\partial_{t}\psi+|\nabla u_{N_{j}}^{q}|^{p-2}\langle\nabla u_{N_{j}}^{q},
		\nabla \psi\rangle}}= 0.$$ By Lemma \ref{weakstrongconv} we have  $$\lim_{j\to \infty}\int_{t_{1}}^{t_{2}}{%
	\int_{\Omega}{u_{N_{j}}\partial_{t}\psi}}=\int_{t_{1}}^{t_{2}}\int_{\Omega}{u\partial_{t}\psi}$$ and $$\lim_{j\to \infty}\left[\int_{\Omega}{u_{N_{j}}\psi}\right]_{t_{1}}^{t_{2}}=\left[\int_{\Omega}{u\psi}\right]_{t_{1}}^{t_{2}},$$ where in the latter we used that $u\in C([0, T], L^{q+1}(\Omega))$ by Lemma \ref{continuitylemmaap}.
Again using the continuity of $u$, we can argue as in Lemma with $\zeta\equiv 1$ and obtain that $\nabla u_{N_{j}}^{q}\to\nabla u^{q}$ in $L^{p}(\Omega\times (0, T))$ as $j\to \infty$. Therefore, using the same arguments as in the proof Lemma \ref{lemexist} we conclude $$\lim_{j\to \infty}\int_{t_{1}}^{t_{2}}\int_{\Omega}|\nabla u_{N_{j}}^{q}|^{p-2}\langle\nabla u_{N_{j}}^{q}, \nabla \psi\rangle=\int_{t_{1}}^{t_{2}}\int_{\Omega}|\nabla u^{q}|^{p-2}\langle\nabla u^{q}, \nabla \psi \rangle,$$ which implies (\ref{usatsifies}).

Testing in (\ref{usatsifies}) with $t_{1}=0$, $t_{2}=\varepsilon$ and test function $\psi(x, t)=\varphi(x)\zeta_{\varepsilon}(t)$, where $\varphi\in C_{0}^{\infty}(\Omega)$ and $$\zeta_{\varepsilon}(t)=\left\{ 
\begin{array}{ll}
	\frac{1}{\varepsilon}(\varepsilon-t)&t\in [0, \varepsilon], \\ 
	0, & t\geq \varepsilon,
\end{array}%
\right. $$ yields $$\frac{1}{\varepsilon}\int_{0}^{\varepsilon}\int_{\Omega}u\varphi\to\int_{\Omega}u(\cdot, 0)\varphi\quad \textnormal{as}~ \varepsilon \to 0.$$ On the other hand, since $u\in C([0, T], L^{q+1}(\Omega))$, we have $$\frac{1}{\varepsilon}\int_{0}^{\varepsilon}\int_{\Omega}u\varphi\to\int_{\Omega}u_{0}\varphi\quad \textnormal{for}~ \varepsilon \to 0$$, which implies the initial condition since $\varphi\in C_{0}^{\infty}(\Omega)$ is arbitrary.
\end{proof}

\section{Estimates of solutions in a cylinder}\label{secmvi}

Let us set $$\delta=q(p-1)-1.$$

\subsection{Caccioppoli type inequalities}

\label{Secini}
\begin{lemma}
	\label{Lem1} Let $u=u\left( x,t\right) $ be a
	non-negative bounded subsolution to \emph{(\ref{Balldtv})} in a cylinder $\Omega\times (0, T)$.
	Let $\eta \left( x,t\right) $ be a locally Lipschitz non-negative bounded
	function in $\Omega\times [0, T]$ such that $\eta \left( \cdot ,t\right) $ has
	compact support in $\Omega $ for all $t\in [0, T]$. Fix some real $\sigma$
	such that 
	\begin{equation}
		\sigma \geq \max(p, pq)  \label{la>3-m}
	\end{equation}%
	and set 
	\begin{equation}
		\lambda =\sigma -\delta \ \ \ \ \text{and}\ \ \ \ \ \alpha =\dfrac{\sigma }{p%
		}.  \label{alpha}
	\end{equation}%
	Choose $0\leq t_{1}<t_{2}\leq T$ and set $Q=\Omega \times \left[ t_{1},t_{2}%
	\right] $. Then, for any $\theta_{1}>\theta_{0}>0$,
	\begin{align}
		\left[ \int_{\Omega }(u-\theta_{1})_{+}^{\lambda }\eta ^{p}\right] _{t_{1}}^{t_{2}}&\nonumber+c_{1}\left(\frac{\theta_{1}}{\theta_{1}-\theta_{0}}\right)^{-(q-1)(p-1)_{-}}
		\int_{Q}\left\vert \nabla \left( (u-\theta_{1})_{+}^{\alpha }\eta \right) \right\vert
		^{p}\\&\leq \int_{Q}\left[ p(u-\theta_{0})_{+}^{\lambda }\eta ^{p-1}\partial _{t}\eta
		+c_{2}\left(\frac{\theta_{1}}{\theta_{1}-\theta_{0}}\right)^{(q-1)(p-1)_{+}}(u-\theta_{0})_{+}^{\sigma }\left\vert \nabla \eta \right\vert ^{p}\right] ,
		\label{veta1}
	\end{align}%
	where $c_{1},c_{2}$ are positive constants depending on $p$, $q$, $\lambda $.
\end{lemma}

Let us recall for later that \begin{equation}\label{valpha}v^{\alpha}=(u-\theta_{1})_{+}^{\alpha}\in L^{p}\left((0, T); W_{0}^{1, p}(\Omega)\right).\end{equation}
Indeed, using $\alpha\geq q$, we get that the function $\Phi(s)=s^{\frac{\alpha}{q}}$ is Lipschitz on any bounded interval in $[0, \infty)$. Thus, $u^{\alpha}=\Phi(u^{q})\in W_{0}^{1, p}(\Omega)$ and $\left|\nabla u^{\alpha}\right|=\left|\Phi^{\prime}(u^{q})\nabla u^{q}\right|\leq C\left|\nabla u^{q}\right|$, whence \begin{equation}\label{valpha0}\int_{Q}u^{\alpha p}+\left\vert \nabla  u^{\alpha } \right\vert^{p}\leq \left(||u||_{L^{\infty}(Q)}^{\sigma-pq}+C\right) \int_{Q}u^{pq}+\left\vert \nabla  u^{q } \right\vert^{p}.\end{equation}
Therefore, $|\nabla v^{\alpha}|=\alpha v^{\alpha-1}|\nabla v|\leq \alpha u^{\alpha-1}=|\nabla u^{\alpha}|$, since $\alpha\geq 1$ by (\ref{la>3-m}) and $v\leq u$.

\begin{proof}
 Since $u$ is a weak subsolution
	of (\ref{Balldtv}), we obtain from Lemma 2.5 in \cite{grigor2023finite}, 
	\begin{equation}  \label{stekinlem}
		\int_{0}^{\tau}\int_{\Omega}{(\partial_{t} u_{h_{0}})\psi+\langle[|\nabla
			u^{q}|^{p-2}\nabla u^{q}]_{h}, \nabla \psi\rangle }\leq 0,
	\end{equation}
	for all $\tau\in (0, T)$ and $\psi\in L^{p}\left((0, \tau);
	W_{0}^{1, p}(\Omega)\right)$.
	It follows from (\ref{la>3-m}) that $\lambda\geq \max(2, 1+q)$. Hence, using the same arguemtents as for (\ref{valpha}), we can test in (\ref{stekinlem}) with the function $\psi=v^{\lambda-1}\eta^{p}$, where  $v=(u-\widetilde{\theta})_{+}$ and $\widetilde{\theta}=(\theta_{1}+\theta_{0})/2$ and can show as in Lemma 2.6 in \cite{grigor2023finite} by using $u\in C([0, T], L^{1}(B))$, that
	\begin{equation}  \label{beforeyoung}
		\left[\int_{\Omega}{v^{\lambda}\eta^{p}}\right]_{t_{1}}^{t_{2}}\leq \int_{Q}{%
			-\lambda\langle|\nabla u^{q}|^{p-2}\nabla u^{q}, \nabla
			(v^{\lambda-1}\eta^{p})\rangle +pv^{\lambda}\eta^{p-1}\partial_{t}\eta}.
	\end{equation}
	We have 
	$\nabla (v^{\lambda -1}\eta ^{p})=(\lambda -1)\eta ^{p}v^{\lambda -2}\nabla
		v+p\eta ^{p-1}v^{\lambda -1}\nabla \eta $.
	Therefore, by (\ref{beforeyoung}), we obtain 
	\begin{align}
		\left[ \int_{\Omega }{v^{\lambda }\eta ^{p}}\right] _{t_{1}}^{t_{2}}&+\int_{Q}\lambda (\lambda -1)q^{p-1}v^{\lambda -2}u^{(q-1)(p-1)}\eta ^{p}|\nabla
		v|^{p}\\ &\leq
		\int_{Q}{\lambda q^{p-1}pv^{\lambda -1}u^{(q-1)(p-1)}|\nabla v|^{p-1}|\nabla \eta |\eta
			^{p-1}} 
		 +\int_{Q}{pv^{\lambda }\eta ^{p-1}\partial _{t}\eta }.  \label{initialfull}
	\end{align}%
	Then by Young's inequality we have, for all $\varepsilon >0$, 
	\begin{align}
		v^{\lambda-1}|\nabla v|^{p-1}|\nabla \eta |\eta ^{p-1}& =\left(
		v^{(\lambda-2)\frac{p-1}{p}}|\nabla v|^{p-1}\eta ^{p-1}\right) \left(
		v^{\frac{\lambda-2+p}{p} }|\nabla \eta |\right)  \notag \\
		& \leq \varepsilon ^{p^{\prime }}v^{\lambda-2}|\nabla v|^{p}\eta ^{p}+%
		\frac{1}{\varepsilon ^{p}}v^{\lambda-2+p}|\nabla \eta |^{p},  \label{youngineq}
	\end{align}%
	where $p^{\prime }=\frac{p}{p-1}$. Combining this with (\ref{initialfull}),
	we deduce 
	\begin{align*}
		\left[ \int_{\Omega }{v^{\lambda }\eta ^{p}}\right] _{t_{1}}^{t_{2}}+&\int_{Q}\lambda q^{p-1}(\lambda -1-p\varepsilon ^{p^{\prime }})u^{(q-1)(p-1)}v^{\lambda-2}|\nabla v|^{p}\eta ^{p}\\&\leq
		\int_{Q}{\frac{\lambda q^{p-1}p}{\varepsilon ^{p}}u^{(q-1)(p-1)}v^{\lambda-2+p}|\nabla \eta |^{p}+pv^{\lambda }\eta ^{p-1}\partial _{t}\eta }.
	\end{align*}%
	In the case $q\geq 1$, we have, since $u\geq v\geq (u-\theta_{1})_{+}$,
	\begin{align*}
		\left\vert \nabla \left( (u-\theta_{1})_{+}^{\alpha }\eta \right) \right\vert ^{p}&\leq 2^{p-1}\alpha ^{p}|\nabla (u-\theta_{1})_{+}|^{p}(u-\theta_{1})_{+}^{p(\alpha -1)}\eta
		^{p}+2^{p-1}(u-\theta_{1})_{+}^{\alpha p}\left\vert \nabla \eta \right\vert ^{p}\\&\leq 2^{p-1}\alpha ^{p}|\nabla v|^{p}u^{(q-1)(p-1)}v^{\lambda-2}\eta
		^{p}+2^{p-1}u^{(q-1)(p-1)}v^{\lambda-2+p}\left\vert \nabla \eta \right\vert ^{p},
	\end{align*}%
	which implies that 
	\begin{equation*}
		|\nabla v|^{p}u^{(q-1)(p-1)}v^{\lambda-2}\eta
		^{p}\geq 2^{1-p}\alpha ^{-p}\left\vert
		\nabla \left( (u-\theta_{1})_{+}^{\alpha }\eta \right) \right\vert ^{p}-\alpha ^{-p}u^{(q-1)(p-1)}v^{\lambda-2+p}\left\vert \nabla \eta \right\vert ^{p}.
	\end{equation*}
	Note that when $u>\theta_{1}>\widetilde{\theta}$ \begin{equation}\label{boundofu}u=u-\widetilde{\theta}+\widetilde{\theta}\leq u-\widetilde{\theta}+\frac{u-\widetilde{\theta}}{\theta_{1}-\widetilde{\theta}}\widetilde{\theta}=\frac{\theta_{0}}{\theta_{0}-\widetilde{\theta}}v=2\frac{\theta_{1}}{\theta_{1}-\theta_{0}}v.\end{equation}
	Therefore, in the case when $q<1$, \begin{align*}\left\vert \nabla \left( (u-\theta_{1})_{+}^{\alpha }\eta \right) \right\vert ^{p}&\leq 2^{p-1}\alpha ^{p}|\nabla v|^{p}v^{p(\alpha -1)}\eta
		^{p}+2^{p-1}v^{\alpha p}\left\vert \nabla \eta \right\vert ^{p}\\&\leq 2^{-q(p-1)}\alpha^{p}|\nabla v|^{p}u^{(q-1)(p-1)}\left(\frac{\theta_{1}}{\theta_{1}-\theta_{0}}\right)^{-(q-1)(p-1)}v^{\lambda-2}\eta
	^{p}\\&\quad+2^{-q(p-1)}u^{(q-1)(p-1)}\left(\frac{\theta_{1}}{\theta_{1}-\theta_{0}}\right)^{-(q-1)(p-1)}v^{\lambda-2+p}\left\vert \nabla \eta \right\vert ^{p}.\end{align*}
	Thus, we obtain in all cases \begin{align}\nonumber|\nabla v|^{p}u^{(q-1)(p-1)}v^{\lambda-2}\eta
		^{p}&\geq \min\left(2^{1-p}, 2^{q(p-1)}\right)\alpha^{-p}\left\vert
		\nabla \left( (u-\theta_{1})_{+}^{\alpha }\eta \right) \right\vert ^{p}\left(\frac{\theta_{1}}{\theta_{1}-\theta_{0}}\right)^{-(q-1)(p-1)_{-}}\\&\quad-\alpha ^{-p}u^{(q-1)(p-1)}v^{\lambda-2+p}\left\vert \nabla \eta \right\vert ^{p}\label{lowerutheta}.\end{align}
	Therefore, 
	\begin{align*}
		\left[ \int_{\Omega }{v^{\lambda }\eta ^{p}}\right] _{t_{1}}^{t_{2}}&+c_{1}\left(\frac{\theta_{1}}{\theta_{1}-\theta_{0}}\right)^{-(q-1)(p-1)_{-}}\int_{Q}{\left\vert \nabla \left( (u-\theta_{1})_{+}^{\alpha }\eta \right) \right\vert ^{p}}\\& \leq c_{2}\int_{Q}{ v^{\lambda-2+p}u^{(q-1)(p-1)}|\nabla
			\eta |^{p}+pv^{\lambda }\eta ^{p-1}\partial _{t}\eta},
	\end{align*}%
	where $c_{1}=\lambda q^{p-1}(\lambda -1-p\varepsilon ^{p^{\prime }})\min\left(2^{1-p}, 2^{q(p-1)}\right)\alpha
		^{-p}$, $c_{2}=\lambda q^{p-1}\left( \left( \lambda -1-p\varepsilon ^{p^{\prime }}\right)
		\alpha ^{-p}+\frac{p}{\varepsilon ^{p}}\right)$.
	Using a similar argument as in (\ref{boundofu}), we have $$ u^{(q-1)(p-1)}\leq 2^{(q-1)(p-1)_{+}}v^{(q-1)(p-1)}\left(\frac{\theta_{1}}{\theta_{1}-\theta_{0}}\right)^{(q-1)(p-1)_{+}}.$$ 
	Hence, we obtain \begin{align}\notag
		\left[ \int_{\Omega }{(u-\theta_{0})_{+}^{\lambda }\eta ^{p}}\right] _{t_{1}}^{t_{2}}&+c_{1}\left(\frac{\theta_{1}}{\theta_{1}-\theta_{0}}\right)^{-(q-1)(p-1)_{-}}\int_{Q}{\left\vert \nabla \left( (u-\theta_{1})_{+}^{\alpha }\eta \right) \right\vert ^{p}}\\&\label{forcorollary} \leq c_{2}^{\prime}\left(\frac{\theta_{1}}{\theta_{1}-\theta_{0}}\right)^{(q-1)(p-1)_{+}}\int_{Q}{ v^{\alpha p}|\nabla
			\eta |^{p}+pv^{\lambda }\eta ^{p-1}\partial _{t}\eta}\\& \notag\leq c_{2}^{\prime}\left(\frac{\theta_{1}}{\theta_{1}-\theta_{0}}\right)^{(q-1)(p-1)_{+}}\int_{Q}{ (u-\theta_{0})_{+}^{\alpha p}|\nabla
			\eta |^{p}+p(u-\theta_{0})_{+}^{\lambda }\eta ^{p-1}\partial _{t}\eta}.
	\end{align}%
	Finally, choosing $\varepsilon $ small enough so that $c_{1}$ and $c_{2}$ are positive we obtain (\ref{veta1}) which finishes the proof.
\end{proof}

\begin{lemma}
	\label{Lem2cacc}Let $u=u\left( x,t\right) $ be a
	non-negative bounded subsolution to \emph{(\ref{Balldtv})} in a cylinder $\Omega\times (0, T)$.
	Let $\sigma$, $\lambda$ and $\alpha$ be as in (\ref{alpha}) and assume that $\sigma\geq pq$.
	Choose $0\leq t_{1}<t_{2}\leq T$ and set $Q=\Omega \times \left[ t_{1},t_{2}%
	\right] $. Then
	\begin{align}
		\left[ \int_{\Omega }u^{\lambda }\right] _{t_{1}}^{t_{2}}+c_{1}
		\int_{Q}\left\vert \nabla  u^{\alpha } \right\vert
		^{p}\leq0,
		\label{veta12}
	\end{align}%
	where $c_{1}$ is a positive constant depending on $p$, $q$, $\lambda $.
\end{lemma}

\begin{proof}
Note that $u\in L^{p}((0, T), W_{0}^{1, p}(\Omega))$ and $u\in C([0, T], L^{1}(\Omega))$. Here, it is sufficient to assume that $\sigma\geq pq$ since this already implies $u^{\alpha}\in L^{p}\left((0, T); W_{0}^{1, p}(\Omega)\right)$ (see (\ref{valpha0})).  Hence, applying the same method as in Lemma 2.6 in \cite{grigor2023finite}, we obtain the claim.
\end{proof}

\begin{lemma}
	\label{monl1}
	Let $u=u\left( x,t\right) $ be a non-negative bounded subsolution to \emph{(\ref{Balldtv})} in $\Omega\times(0, T)$. If $\lambda\geq 1$, including $\lambda=\infty$, then the function 
	\begin{equation*}
		t\mapsto \left\Vert u(\cdot ,t)\right\Vert _{L^{\lambda}(\Omega)}
	\end{equation*}%
	is monotone decreasing in $[0, T]$.
\end{lemma}

\begin{proof}
We have $u\in L^{p}((0, T), W_{0}^{1, p}(\Omega))$ and $u\in C([0, T], L^{1}(\Omega))$. Hence, applying the same method as in Lemma 2.9 in \cite{Grigoryan2024}, the claim follows.
\end{proof}

\subsection{Comparison in two cylinders}

\begin{lemma}
	\label{Lem2yyy}Consider two precompact balls $B_{0}=B\left( x_{0},r_{0}\right) $ and $%
	B_{1}=B\left( x_{0},r_{1}\right) $ with $0<r_{1}<r_{0}$ where $B_{0}$ is
	precompact. Consider two cylinders $%
	Q_{i}=B_{i}\times \lbrack 0,T]$, $i=0,1.$ 
	Let $u$ be a non-negative bounded subsolution in $Q_{0}$ such that \begin{equation}\label{u000}
		u\left( \cdot ,0\right) =0.\end{equation} For $\theta_{1}>\theta_{0}>0$, denote
	\begin{equation*}
		v_{i}=\left( u-\theta_{i} \right) _{+}.
	\end{equation*}%
	Let $\sigma $ and $\lambda $ be reals satisfying (\ref{la>3-m}) and (\ref{alpha}).
	Set 
	\begin{equation*}
		J_{i}=\int_{Q_{i}}v_{i}^{\sigma }d\mu dt.
	\end{equation*}%
	Then%
	\begin{equation}\label{compsingsp1}
		J_{1}\leq \frac{Cr_{0}^{p}}{\left( \iota (B_{0})\mu (B_{0})\left( r_{0}-r_{1}\right) ^{p}\right) ^{\nu
			}(\theta_{1}-\theta_{0}) ^{\lambda \nu }\left( r_{0}-r_{1}\right) ^{p}}\left(\frac{\theta_{1}}{\theta_{1}-\theta_{0}}\right)^{|q-1|(p-1)+\nu(q-1)(p-1)_{+}}J_{0} ^{1+\nu },
	\end{equation}%
	where $\nu $ is the Faber-Krahn exponent, $\iota (B_{0})$ is the Faber-Krahn
	constant in $B_{0}$, and $C$ depends on $p$, $q$ and $\lambda $.
	\end{lemma}

\begin{proof}
	Let $\eta(x, t)=\eta \left( x\right) $ be a bump function of $B_{1}$ in $B_{1/2}:=B\left( x_{0},\frac{r_{0}+r_{1}}{2}\right)$. Recall that by (\ref{valpha}), $v_{1}^{\alpha }\eta\in L^{p}\left([0,T]; W_{0}^{1, p}(B)\right)$, where $\alpha$ is defined by (\ref{alpha}), that is $\alpha=\frac{\sigma}{p}$. Hence, applying the Faber-Krahn inequality in ball $B_{0}$ for any $t\in [0, T]$ we get that
	\begin{equation}\label{FKapp2sing}
		\int_{B_{1}}v_{1}^{\sigma}\leq \int_{B_{0}}\left( v_{1}^{\alpha }\eta \right) ^{p}\leq r_{0}^{p}\left( \frac{\mu
			\left( D_{t}\right) }{\iota (B_{0})\mu (B_{0})}\right) ^{\nu
		} \int_{B_{0}}\left\vert \nabla \left( v_{1}^{\alpha }\eta \right) \right\vert^{p},
	\end{equation}%
	where $\iota(B_{0})$ is the \textit{normalized Sobolev constant} in $B_{0}$ and we used that $\alpha p=\sigma$ and $\eta =1$ in $B_{1}$ and \begin{equation*}
		D_{t}=\left\{ v_{1}^{\alpha }\eta \left( \cdot ,t\right) >0\right\} =\left\{
		v_{1}>0\right\} \cap \left\{ \eta >0\right\} =\left\{ u\left( \cdot,t\right)
		>\theta_{1} \right\} \cap B_{1/2}.
	\end{equation*} Also, note that $\eta_{t}=0$ and $\left\vert \nabla \eta \right\vert \leq \frac{2}{r_{0}-r_{1}}$. From (\ref{veta1}) we therefore obtain 
	\begin{align}\notag
		c_{1}\int_{0}^{T}\int_{B_{0}}\left\vert \nabla \left( v_{1}^{\alpha }\eta
		\right) \right\vert ^{p}&\leq c_{2}\left(\frac{\theta_{1}}{\theta_{1}-\theta_{0}}\right)^{|q-1|(p-1)}
		\int_{0}^{T}\int_{B_{0}}v_{0}^{\sigma}|\nabla \eta|^{p}\\&\label{gradalpha3}\leq \frac{c_{3}}{\left(r_{0}-r_{1}\right) ^{p}}\left(\frac{\theta_{1}}{\theta_{1}-\theta_{0}}\right)^{|q-1|(p-1)}J_{0},
	\end{align} where $c_{3}=c_{2}2^{p}$.
	
	Let us now apply Lemma \ref{Lem1} to function $v_{0}$ in $B_{0}\times \left[0, T\right] $. Take $\eta \left( x,t\right) =\eta \left( x\right)$,
	where $\eta$ is a bump function of $B_{1/2}$ in $B_{0}$ so that $\left\vert \nabla \eta \right\vert \leq \frac{2}{r_{0}-r_{1}}$.
	From (\ref{veta1}) we obtain
	\begin{align*}
		\left[\int_{B_{0}}v_{1}^{\lambda}\eta ^{p}\right] _{0}^{T}\leq 
		\int_{0}^{T}\int_{B_{0}} c_{2}\left(\frac{\theta_{1}}{\theta_{1}-\theta_{0}}\right)^{(q-1)(p-1)_{+}}\left\vert \nabla \eta \right\vert ^{p}v_{0}^{\sigma } .
	\end{align*}%
	Hence, by (\ref{u000}) and since $\eta(x)=1$ for $x\in B_{1/2}$, $$\int_{B_{1/2}}v_{1}^{\lambda }\left( \cdot ,t\right)\leq \frac{c_{3}}{(r_{0}-r_{1})^{p}}\left(\frac{\theta_{1}}{\theta_{1}-\theta_{0}}\right)^{(q-1)(p-1)_{+}}J_{0}.$$ 
	Thus, we deduce
	\begin{equation*}
		\mu \left( D_{t}\right) \leq \frac{1}{(\theta_{1}-\theta_{0})^{\lambda }}%
		\int_{B_{1/2}\cap\{u>\theta_{1}\}}v_{1}^{\lambda }\left( \cdot ,t\right) \leq \frac{c_{3}}{(r_{0}-r_{1})^{p}
			(\theta_{1}-\theta_{0}) ^{\lambda}}\left(\frac{\theta_{1}}{\theta_{1}-\theta_{0}}\right)^{(q-1)(p-1)_{+}}J_{0}.
	\end{equation*}%
	
	Combining this with (\ref{FKapp2sing}) and (\ref{gradalpha3}) we obtain
	\begin{align*}
		J_{1}=\int_{t_{1}}^{T}\int_{B_{1}}v_{1}^{\sigma}&\leq r_{0}^{p}\left( \frac{c_{3}J_{0}}{\iota (B_{0})\mu (B_{0})(r_{0}-r_{1})^{p}(\theta_{1}-\theta_{0}) ^{\lambda }}\left(\frac{\theta_{1}}{\theta_{1}-\theta_{0}}\right)^{(q-1)(p-1)_{+}}\right) ^{\nu}\\&\quad\times
		\frac{c_{3}}{c_{1}\left(r_{0}-r_{1}\right) ^{p}}\left(\frac{\theta_{1}}{\theta_{1}-\theta_{0}}\right)^{|q-1|(p-1)}J_{0},
	\end{align*}%
	which implies (\ref{compsingsp1}) and finishes the proof.
\end{proof}

\subsection{Iterations and the mean value theorem}

\begin{lemma}
	\label{Tmeansing}Let the ball $B=B\left( x_{0},R\right) $ be precompact. Let $u$
	be a non-negative bounded subsolution in 
	$Q=B\times \left[ 0,t\right]$ such that \begin{equation*}
		u\left( \cdot ,0\right) =0.\end{equation*}
	Let $\sigma$ and $\lambda$ be reals such that \begin{equation}\label{silapossing}\sigma>0\quad \textnormal{and}\quad \lambda=\sigma-\delta>0.\end{equation}
	Then, for the cylinder 
	$Q^{\prime }=\frac{1}{2}B\times [ 0,t] ,$
	we have%
	\begin{equation}
		\left\Vert u\right\Vert _{L^{\infty }\left( Q^{\prime }\right) }\leq \left( 
		\frac{C}{\iota (B)\mu (B)R^{p}}
			\int_{Q}u^{\sigma}
		\right) ^{1/\lambda } ,  \label{meansing}
	\end{equation}%
	where $\iota (B)$ is the Faber-Krahn constant in $B$,
	and the constant $C$ depends on $p$, $q$, $\lambda $ and $\nu$.

\end{lemma}
\FRAME{ftbpF}{2.1392in}{1.6447in}{0pt}{\Qcb{Cylinders $Q$ and $Q^{\prime }$}}{\Qlb{pic3}}{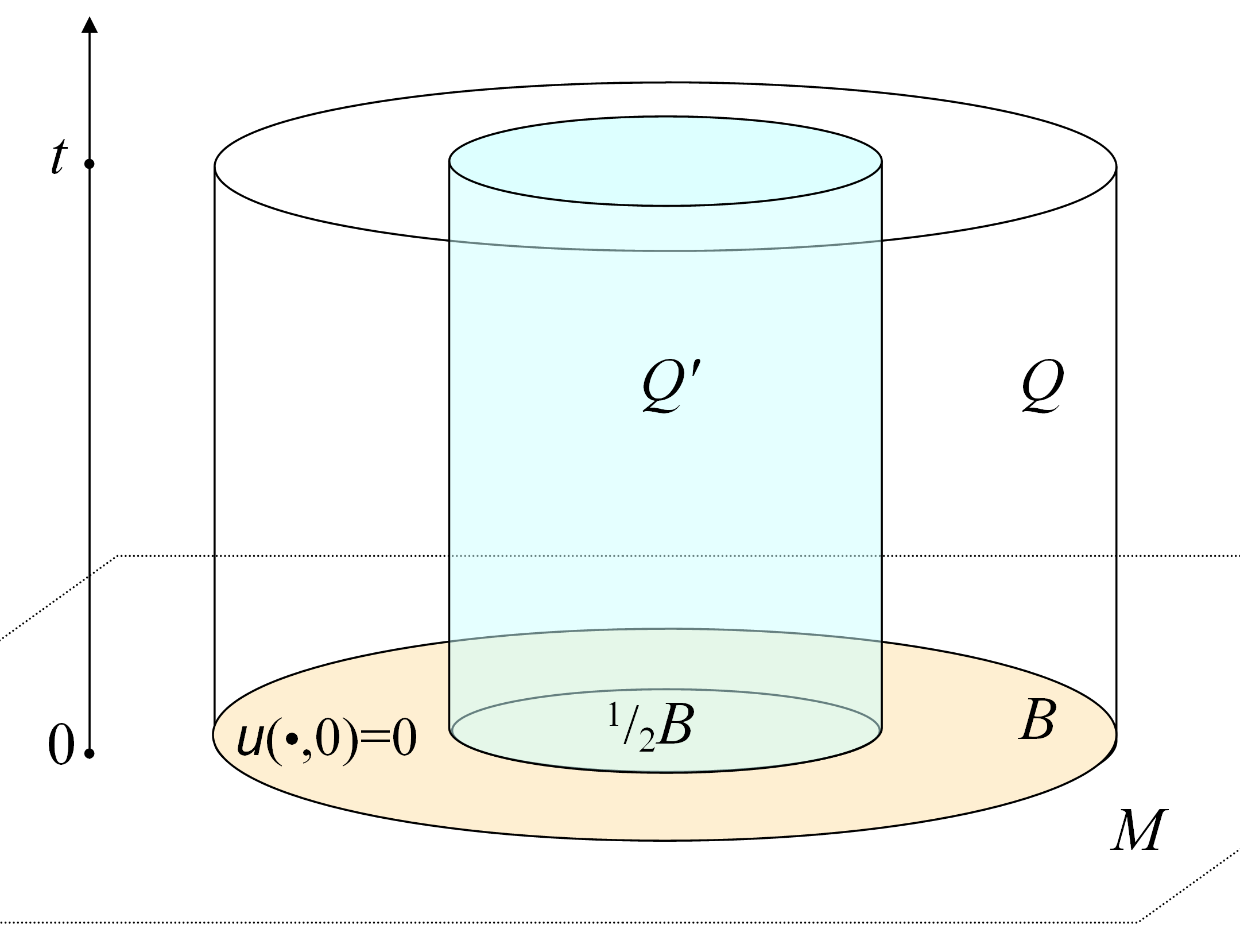}{\special%
	{language "Scientific Word";type "GRAPHIC";maintain-aspect-ratio
		TRUE;display "USEDEF";valid_file "F";width 2.1392in;height 1.6447in;depth
		0pt;original-width 10.0838in;original-height 7.7406in;cropleft "0";croptop
		"1";cropright "1";cropbottom "0";filename 'pic3.png';file-properties
		"XNPEU";}}
	
\begin{remark}
In \cite{Grigoryan2024} the same mean value inequality was proved under the condition that \begin{equation}\label{allcases}p>2\quad\textnormal{and}\quad\frac{1}{p-1}<q\leq1\quad \textnormal{or}\quad 1<p<2\quad\textnormal{and}\quad1\leq q< \frac{1}{p-1}.\end{equation}
\end{remark}

\begin{proof}
	Let us first prove (\ref{meansing}) for $\sigma$ large enough as in Lemma \ref{Lem1}. Consider sequences 
	\begin{equation*}
		r_{k}=\left( \frac{1}{2}+2^{-k-1}\right) R,
	\end{equation*}%
	where $k=0,1,2,...$, so that
	$r_{0} =R$ and $r_{k}\searrow \frac{1}{2}R$ as $k\rightarrow
	\infty $ .
	Set 
	$B_{k}=B\left( x_{0},r_{k}\right)$, $Q_{k}=B_{k}\times \left[ 0,t   \right]$
	so that 
	$B_{0}=B$, $Q_{0}=Q$ and $Q_{\infty}:=\lim_{k\rightarrow \infty }Q_{k}=Q^{\prime }$.
	        \FRAME{dtbpF}{2.5604in}{1.6587in}{0pt}{\Qcb{Cylinders $Q_{k}$}}{}{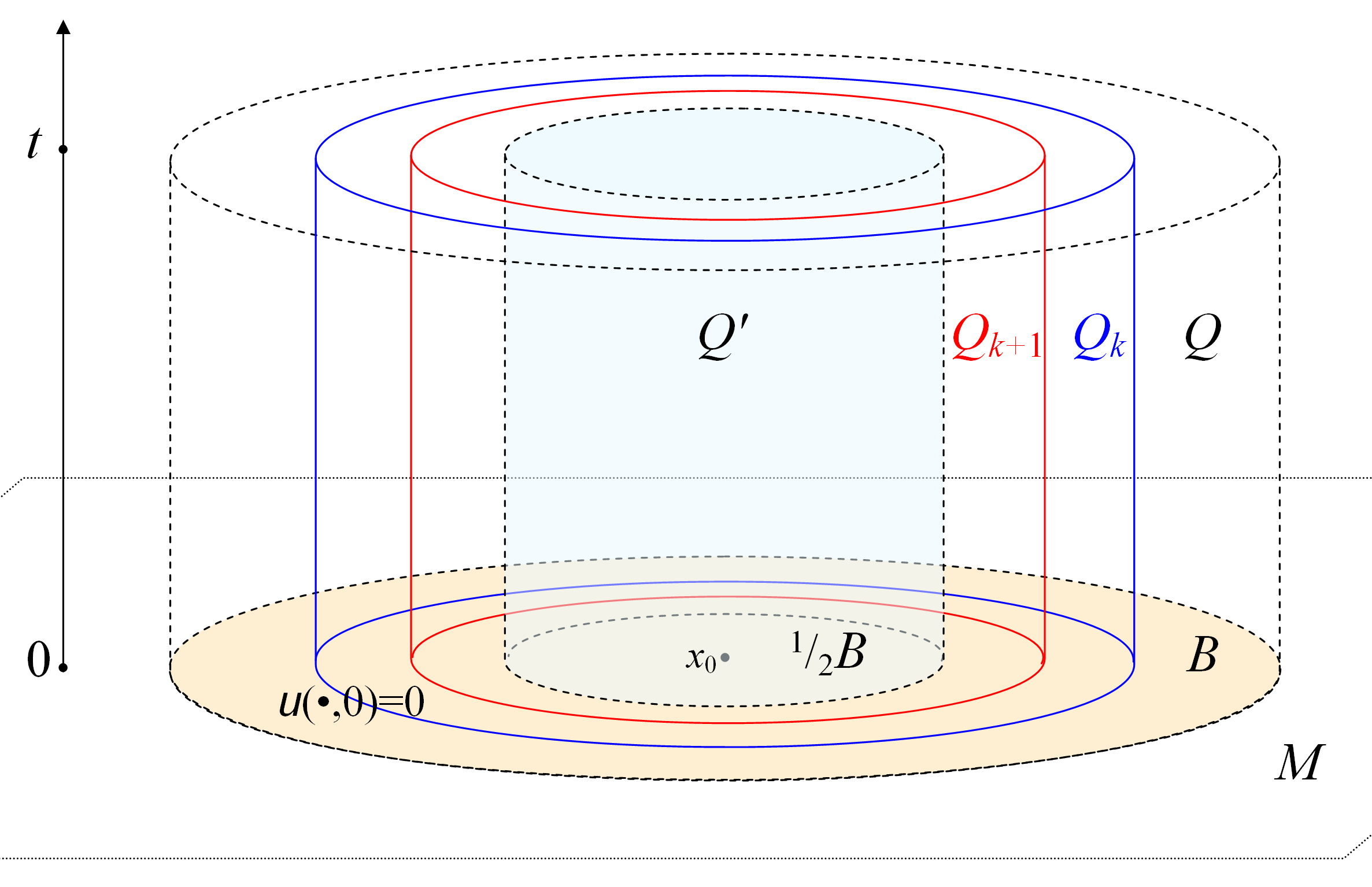}{\special{language
			"Scientific Word";type "GRAPHIC";maintain-aspect-ratio TRUE;display
			"USEDEF";valid_file "F";width 2.5604in;height 1.6587in;depth
			0pt;original-width 11.078in;original-height 7.1615in;cropleft "0";croptop
			"1";cropright "1";cropbottom "0";filename 'pic5.png';file-properties
			"XNPEU";}}
	Choose some $\theta >0$ to be specified later and define $\theta_{k}=\left( 1-2^{-k}\right) \theta$ and
	\begin{equation*}
		u_{k}=(u-\theta_{k})_{+}=\left( u-\left( 1-2^{-k}\right) \theta \right) _{+}.
	\end{equation*}
	Set also $J_{k}=\int_{Q_{k}}u_{k}^{\sigma }$. Then we obtain by Lemma \ref{Lem2yyy} that%
	\begin{align*}
		J_{k+1}&\leq \frac{Cr_{k}^{p}}{\left( \iota (B_{k})\mu (B_{k})\left( r_{k}-r_{k+1}\right) ^{p}\right) ^{\nu
			}(\theta_{k+1}-\theta_{k}) ^{\lambda \nu }\left( r_{k}-r_{k+1}\right) ^{p}}\\&\quad \times\left(\frac{\theta_{k+1}}{\theta_{k+1}-\theta_{k}}\right)^{|q-1|(p-1)+\nu(q-1)(p-1)_{+}}J_{k} ^{1+\nu }.
	\end{align*}%
	By monotonicity we have $\frac{r_{k}^{p}}{\left( \iota (B_{k})\mu (B_{k})\right) ^{\nu }}\leq \frac{%
			R^{p}}{\left( \iota (B)\mu (B)\right) ^{\nu }}$.
	Since $r_{k}-r_{k+1}=2^{-k-2}R$ and $\theta_{k+1}-\theta_{k}=2^{-k-1}\theta$ we get
	\begin{eqnarray*}
		J_{k+1} &\leq &\frac{C2^{k\left(\lambda \nu + p\left( 1+\nu \right)+ |q-1|(p-1)+\nu(q-1)(p-1)_{+} \right) }}{\left( \iota (B)\mu
			(B)R^{p}\right) ^{\nu }\theta ^{\lambda \nu }}J_{k}^{1+\nu }=\frac{A^{k}J_{k}^{1+\nu }}{\Theta }
	\end{eqnarray*}%
	where $A=2^{\left(\lambda \nu + p\left( 1+\nu \right)+ |q-1|(p-1)+\nu(q-1)(p-1)_{+} \right)}\geq 1$ and $\Theta =c\left(\iota (B)\mu (B)\theta ^{\lambda }R^{p}\right)^{\nu}$.
	
	Now let us apply Lemma 6.1 from \cite{grigor2023finite} with $\omega =\nu $: if 
	\begin{equation}
		\Theta \geq A^{1/\nu }J_{0}^{\nu },  \label{Thetayyy}
	\end{equation}%
	then, for all $k\geq 0,$%
	\begin{equation*}
		J_{k}\leq A^{-k/\nu }J_{0}.
	\end{equation*}%
	In terms of $\theta $ the condition (\ref{Thetayyy}) is equivalent%
	\begin{equation*}
		c\left(\iota (B)\mu (B)\theta ^{\lambda }R^{p}\right)^{\nu}\geq A^{1/\nu }J_{0}^{\nu }
	\end{equation*}%
	that is,%
	\begin{equation*}
		\theta \geq \left(\frac{CJ_{0}}{\iota (B)\mu (B)R^{p}}\right)^{1/\lambda}.
	\end{equation*}%
	Hence, we choose $\theta $ to have equality here and for this $\theta $ we
	obtain $J_{k}\rightarrow 0$ as $k\rightarrow \infty $, which implies that
	$u\leq \theta$ in $Q_{\infty }.$
	Hence,%
	\begin{equation}\label{uppermitStheta}
		\left\Vert u\right\Vert _{L^{\infty }\left( Q^{\prime }\right) }\leq \left( 
		\frac{CJ_{0}}{\iota (B)\mu (B)R^{p}}%
		\right) ^{1/\lambda },
	\end{equation}
	which proves (\ref{meansing}) for any large enough $\sigma$.
	
	By standard methods (see \cite{Grigoryan2024}) we conclude from this that (\ref{meansing}) holds for any $\sigma >0$.	
	\end{proof}

\section{Global existence theorem}\label{secglobex}

The next result contains Theorem \ref{mainthmint} from the Introduction.

\begin{theorem}\label{mainexthm}
Assume that $u_{0}\in L^{1}(M)\cap L^{\infty}(M)$ is non-negative. Then there exists a unique non-negative bounded weak solution $u$ of the Cauchy problem \begin{equation}
	\left\{ 
	\begin{array}{ll}
		\partial _{t}u=\Delta _{p}u^{q} &\text{in}~M\times (0, \infty), \\u(x, 0)=u_{0}(x)& \text{in}~M,
	\end{array}%
	\right.  \label{wholemandtv}
\end{equation} in the sense of Definition \ref{defweaksolu}.
\end{theorem}

\begin{proof}
Let $B$ be a ball in $M$ and $u_{B}$ be a non-negative solution to the Cauchy problem (\ref{Balldtv}) in $Q=B\times(0, T)$. 
Applying the Caccioppoli type inequality (\ref{veta12}) with $\lambda=1+q$ we get $$\left[ \int_{B }u_{B}^{1+q }\right] _{0}^{T}+c_{1}
\int_{Q}\left\vert \nabla  u_{B}^{q } \right\vert
^{p}\leq 0,$$ and therefore $$\int_{Q}|\nabla u_{B}^{q}|^{p}\leq \frac{1}{c_{1}} \int_{M }u_{0}^{1+q }.$$

Hence, after extending $u_{B}$ to $M$ by setting $u\equiv0$ outside $B$, $u_{B}^{q}$ is a bounded sequence in $L^{p}((0, T), D^{1, p}(M))$. Thus, there is a subsequence $u_{B_{k}}^{q}$ and a function $w\in L^{p}((0, T), D^{1, p}(M))$ such that $u_{B_{k}}^{q}\to w$ weakly in $L^{p}((0, T), D^{1, p}(M))$ as $B_{k}\to M$.
Since the sequence $\{u_{B}\}$ is monotone increasing in $B$ by (\ref{forcomparison}), the limit is also a pointwise limit as $B_{k}\to M$. From Lemma \ref{monl1} we know that, for any $\lambda\geq1$, any ball $B$ in $M$ and any $0< t< T$, \begin{equation}\label{monoinex}||u_{B}(\cdot, t)||_{L^{\lambda}(B)}\leq ||u_{0}||_{L^{\lambda}(M)}.\end{equation} Thus, we obtain from (\ref{monoinex}) by sending $\lambda\to \infty$ that $u:=w^{1/q}$ is bounded. Using (\ref{monoinex}) with $\lambda=q+1$ and Fatou's Lemma, we also obtain that $u\in L^{q+1}(M\times (0, T))$. Therefore, we can argue as in Lemma \ref{continuitylemmaap} and get that $u\in C([0, T], L^{q+1}(M))$.

Using similar arguments as in the proof of Lemma \ref{lemmaexball}, we therefore obtain that $u$ also satisfies (\ref{defvonweaksolq}), that is, for any $0\leq t_{1}<t_{2}\leq T$,  \begin{equation}\label{solinmainpr}\left[\int_{M}{u\psi}\right]_{t_{1}}^{t_{2}}+\int_{t_{1}}^{t_{2}}{%
		\int_{M}{-u\partial_{t}\psi+|\nabla u^{q}|^{p-2}\langle\nabla u^{q},
			\nabla \psi\rangle}}= 0,\end{equation}
for any $\psi\in W^{1, 1+\frac{1}{q}}\left((0, T);L^{1+\frac{1}{q}}(M)\right)
\cap L^{p}\left((0, T); D^{1, p}(M)\right)$ and that $u(\cdot, 0)=u_{0}$.

From Lemma \ref{thmcomp} we get a unique solution in $M\times (0, \infty)$, which finishes the proof.
\end{proof}

\section{Finite propagation speed}\label{secfps}

\subsection{Propagation inside a ball}

The next result contains Theorem \ref{mainthmintprop} from the Introduction.

Using the mean value inequality \ref{Tmeansing} and the method from \cite{Grigoryan2024} we can prove the following result, which improves the range of $p$ and $q$ obtained in Theorem 1.1 from \cite{Grigoryan2024}.

\begin{theorem}
	\label{TFPSq} Assume that $\delta=q(p-1)-1>0$. Let $u$ be a bounded non-negative subsolution in $M\times \left[
	0,T\right] $. Let $B_{0}=B\left( x_{0},R\right) $ be a
	ball such that 
	\begin{equation*}
		u_{0}=0\ \text{in\ }B_{0}.
	\end{equation*}%
	Let $\sigma $ be a real such that 
	\begin{equation}
		\sigma \geq 1\text{ and }\sigma >\delta.  \label{si>>}
	\end{equation}%
	Set%
	\begin{equation}
		t_{0}=\eta \iota (B_{0})\mu (B_{0})^{\frac{\delta}{\sigma }}R^{p}
		||u_{0}||_{L^{\sigma }(M)} ^{-\delta}\wedge T,
		\label{t0q}
	\end{equation}%
	where $\eta =\eta (p,q,\nu ,\sigma )>0$ is sufficiently small. Then%
	\begin{equation*}
		u=0\ \ \text{in\ \ }\frac{1}{2}B_{0}\times \left[ 0,t_{0}\right] .
	\end{equation*}
\end{theorem}

\begin{remark}
Note that in the result in \cite{Grigoryan2024} we also, besides $\delta>0$, assumed that (see (\ref{allcases})) $$p>2\quad\textnormal{and}\quad q\leq 1.$$  
\end{remark}

\subsection{Curvature and propagation rate}
For any $K\subset M$ and any $r>0$, denote by $K_{r}$ a closed $r$-neighborhood of $K.$

Assume that $\textnormal{supp}~u_{0}$ is compact. Then the function $$\rho:\left( 0,\infty\right) \rightarrow \mathbb{R}_{+}$$ is called a \emph{propagation rate} or \emph{propagation function} of $u$ if $\limfunc{supp}%
u\left( \cdot ,t\right) \subset (\textnormal{supp}~u_{0})_{\rho \left( t\right) }\ $for all $t\in
\left( 0,\infty\right) .$

From Theorem \ref{TFPSq} we obtain as in \cite{Grigoryan2024} the following result.

\begin{corollary}\label{Corpropric}
	Assume that $\delta>0$. Let $M$ satisfy the relative Faber-Krahn inequality. Fix a reference point $x_{0}\in \textnormal{supp}~u_{0}$. Suppose that
	for some $\alpha>0$ and all large enough $r$,%
	\begin{equation}
		\mu \left( B\left( x_{0},r\right) \right) \geq cr^{\alpha }.  \label{ra}
	\end{equation}%
	Then $u$ has a propagation function 
	\begin{equation*}
		\rho (t)=Ct^{1/\beta }
	\end{equation*}%
	for large $t$, where
	\begin{equation*}
		\beta =p+\alpha \frac{\delta}{\sigma }
	\end{equation*}%
	with $\sigma $ as in \emph{(\ref{si>>})} and $C$ depends on $\left\Vert u_{0}\right\Vert _{L^{\sigma }(M)},p,q,n,\alpha$ and $c$.
\end{corollary}

\begin{remark}
	In $\mathbb{R}^{n}$ we have (\ref{ra}) with $\alpha=n$. If $\sigma=1$, we obtain the sharp propagation rate $1/\beta$, where $\beta=p+n\delta$. By (\ref{si>>}), we can take $\sigma =1$ provided $\delta<1,$ that is, when $q<\frac{2}{p-1}$. Hence, in the range \begin{equation}\label{sharprate} p>1\quad \textnormal{and}\quad\frac{1}{p-1}<q< \frac{2}{p-1}\end{equation}
	(see Fig. \ref{pic9}), we get a sharp propagation rate, which improves the range $$p>2\quad \textnormal{and}\quad \frac{1}{p-1}<q< \min\left(\frac{2}{p-1}, 1\right),$$ which was obtained in Corollary 4.4 in \cite{Grigoryan2024}. \FRAME{ftbpFU}{3.409in}{1.4878in}{0pt}{\Qcb{Range of $p, q$}}{\Qlb{%
			pic9}}{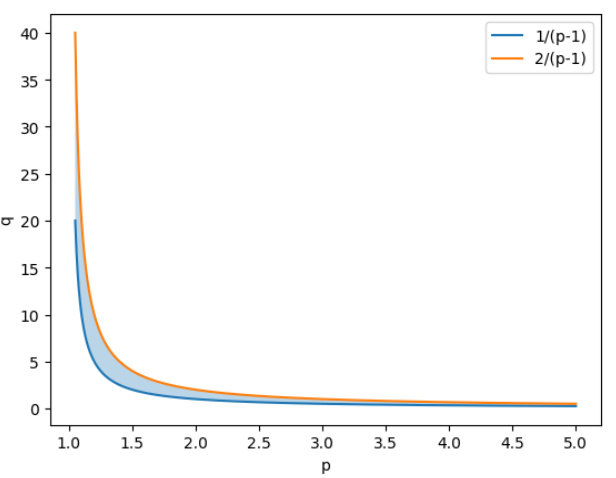}{\special{language "Scientific Word";type
			"GRAPHIC";maintain-aspect-ratio TRUE;display "USEDEF";valid_file "F";width
			3.409in;height 1.4878in;depth 0pt;original-width 13.389in;original-height
			5.8168in;cropleft "0";croptop "1";cropright "1";cropbottom "0";filename
			'pic9.png';file-properties "XNPEU";}} 
	\end{remark}
		\bibliographystyle{abbrv}
	\bibliography{gradleib}

\begin{thebibliography}{10}

\bibitem{andreucci2015optimal}
D.~Andreucci, A.~F. Tedeev, et~al.
\newblock Optimal decay rate for degenerate parabolic equations on noncompact
  manifolds.
\newblock {\em Methods Appl. Anal}, 22(4):359--376, 2015.

\bibitem{anttila2026uniqueness}
R.~Anttila, P.~Lindqvist, and M.~Parviainen.
\newblock Uniqueness of sign-changing solutions to {T}rudinger's equation.
\newblock {\em arXiv preprint arXiv:2602.04748}, 2026.

\bibitem{barbu2025leibenson}
V.~Barbu, S.~Grube, M.~Rehmeier, and M.~R{\"o}ckner.
\newblock The {L}eibenson process.
\newblock {\em arXiv preprint arXiv:2508.12979}, 2025.

\bibitem{bogelein2021existence}
V.~B{\"o}gelein, N.~Dietrich, and M.~Vestberg.
\newblock Existence of solutions to a diffusive shallow medium equation.
\newblock {\em Journal of Evolution Equations}, 21(1):845--889, 2021.

\bibitem{bogelein2018higher}
V.~B{\"o}gelein, F.~Duzaar, R.~Korte, and C.~Scheven.
\newblock The higher integrability of weak solutions of porous medium systems.
\newblock {\em Advances in Nonlinear Analysis}, 8(1):1004--1034, 2018.

\bibitem{bogelein2013parabolic}
V.~B{\"o}gelein, F.~Duzaar, and P.~Marcellini.
\newblock Parabolic systems with $p, q$-growth: a variational approach.
\newblock {\em Archive for Rational Mechanics and Analysis}, 210(1):219--267,
  2013.

\bibitem{bonforte2008fast}
M.~Bonforte, G.~Grillo, and J.~L. Vazquez.
\newblock Fast diffusion flow on manifolds of nonpositive curvature.
\newblock {\em Journal of Evolution Equations}, 8:99--128, 2008.

\bibitem{Buser}
P.~Buser.
\newblock A note on the isoperimetric constant.
\newblock {\em Ann. Sci. Ecole Norm. Sup.}, 15:213--230, 1982.

\bibitem{de1957sulla}
E.~De~Giorgi.
\newblock {S}ulla differenziabilit{\`a} e l'analiticit{\`a} delle estremali
  degli integrali multipli regolari.
\newblock {\em Mem. Accad. Sci. Torino}, 3:25--43, 1957.

\bibitem{de2022wasserstein}
N.~De~Ponti, M.~Muratori, and C.~Orrieri.
\newblock Wasserstein stability of porous medium-type equations on manifolds
  with {R}icci curvature bounded below.
\newblock {\em Journal of Functional Analysis}, 283(9):109661, 2022.

\bibitem{dekkers2005finite}
S.~Dekkers.
\newblock Finite propagation speed for solutions of the parabolic $p$-laplace
  equation on manifolds.
\newblock {\em Communications in Analysis and Geometry}, 13(4):741--768, 2005.

\bibitem{grigor}
A.~Grigor'yan.
\newblock The heat equation on non-compact {R}iemannian manifolds.
\newblock {\em Math. USSR Sb.}, 72:47--77, 1992.

\bibitem{grigor2023finite}
A.~Grigor'yan and P.~S{\"u}rig.
\newblock Finite propagation speed for {L}eibenson’s equation on {R}iemannian
  manifolds.
\newblock {\em Comm. Anal. Geom.}, 2024.

\bibitem{Grigoryan2024}
A.~Grigor'yan and P.~S{\"u}rig.
\newblock Sharp propagation rate for {L}eibenson’s equation on {R}iemannian
  manifolds.
\newblock {\em Ann. Scuola Norm. Super. Pisa}, 2024.

\bibitem{grigor2025upper}
A.~Grigor'yan and P.~S{\"u}rig.
\newblock Upper bounds for solutions of {L}eibenson's equation on {R}iemannian
  manifolds.
\newblock {\em Journal of Functional Analysis}, page 110878, 2025.

\bibitem{grillo2025porous}
G.~Grillo, D.~D. Monticelli, and F.~Punzo.
\newblock The porous medium equation on noncompact manifolds with nonnegative
  {R}icci curvature: {A} {G}reen function approach.
\newblock {\em Journal of Differential Equations}, 430:113191, 2025.

\bibitem{grillo2016smoothing}
G.~Grillo and M.~Muratori.
\newblock Smoothing effects for the porous medium equation on
  {C}artan--{H}adamard manifolds.
\newblock {\em Nonlinear Analysis}, 131:346--362, 2016.

\bibitem{grillo2018porous}
G.~Grillo, M.~Muratori, and F.~Punzo.
\newblock The porous medium equation with measure data on negatively curved
  {R}iemannian manifolds.
\newblock {\em Journal of the European Mathematical Society},
  20(11):2769--2812, 2018.

\bibitem{grillo2021fast}
G.~Grillo, M.~Muratori, and F.~Punzo.
\newblock Fast diffusion on noncompact manifolds: well-posedness theory and
  connections with semilinear elliptic equations.
\newblock {\em Transactions of the American Mathematical Society},
  374(9):6367--6396, 2021.

\bibitem{grillo2017porous}
G.~Grillo, M.~Muratori, and J.~L. V{\'a}zquez.
\newblock The porous medium equation on {R}iemannian manifolds with negative
  curvature. the large-time behaviour.
\newblock {\em Advances in Mathematics}, 314:328--377, 2017.

\bibitem{ishige1996existence}
K.~Ishige.
\newblock On the existence of solutions of the cauchy problem for a doubly
  nonlinear parabolic equation.
\newblock {\em SIAM Journal on Mathematical Analysis}, 27(5):1235--1260, 1996.

\bibitem{ivanov1992existence}
A.~Ivanov and P.~Mkrtychyan.
\newblock Existence of {H}{\"o}lder continuous generalized solutions of the
  first boundary value problem for quasilinear doubly degenerate parabolic
  equations.
\newblock {\em Journal of Soviet Mathematics}, 62(3):2725--2740, 1992.

\bibitem{ivanov1997regularity}
A.~V. Ivanov.
\newblock Regularity for doubly nonlinear parabolic equations.
\newblock {\em Journal of Mathematical Sciences}, 83(1):22--37, 1997.

\bibitem{kinnunen2006pointwise}
J.~Kinnunen and P.~Lindqvist.
\newblock Pointwise behaviour of semicontinuous supersolutions to a quasilinear
  parabolic equation.
\newblock {\em Annali di Matematica Pura ed Applicata}, 185(3):411--435, 2006.

\bibitem{ladyzhenskaya1968linear}
O.~Ladyzhenskaya, V.~Solonnikov, and N.~Ural’tseva.
\newblock Linear and quasilinear equations of parabolic type, transl. math.
\newblock {\em Monographs, Amer. Math. Soc}, 23, 1968.

\bibitem{leibenzon1945general}
L.~Leibenson.
\newblock General problem of the movement of a compressible fluid in a porous
  medium. izv akad. nauk sssr.
\newblock {\em Geography and Geophysics}, 9:7--10, 1945.

\bibitem{leibenson1945turbulent}
L.~Leibenson.
\newblock Turbulent movement of gas in a porous medium.
\newblock {\em Izv. Akad. Nauk SSSR Ser. Geograf. Geofiz}, 9:3--6, 1945.

\bibitem{meglioli2025global}
G.~Meglioli, F.~Oliva, and F.~Petitta.
\newblock Global existence for a {L}eibenson type equation with reaction on
  {R}iemannian manifolds.
\newblock {\em arXiv preprint arXiv:2505.08304}, 2025.

\bibitem{otto1996l1}
F.~Otto.
\newblock L1-contraction and uniqueness for quasilinear elliptic--parabolic
  equations.
\newblock {\em Journal of differential equations}, 131(1):20--38, 1996.

\bibitem{raviart1970resolution}
P.-A. Raviart.
\newblock Sur la r{\'e}solution de certaines {\'e}quations paraboliques non
  lin{\'e}aires.
\newblock {\em Journal of Functional Analysis}, 5(2):299--328, 1970.

\bibitem{Saloff}
L.~Saloff-Coste.
\newblock {\em Aspects of Sobolev-type inequalities}.
\newblock LMS Lecture Notes Series, vol. 289. Cambridge Univ. Press, 2002.

\bibitem{surig2024finite}
P.~S{\"u}rig.
\newblock Finite extinction time for subsolutions of the weighted {L}eibenson
  equation on {R}iemannian manifolds.
\newblock {\em arXiv preprint arXiv:2412.06496}, 2024.

\bibitem{surig2024sharp}
P.~S{\"u}rig.
\newblock Sharp sub-{G}aussian upper bounds for subsolutions of {T}rudinger’s
  equation on {R}iemannian manifolds.
\newblock {\em Nonlinear Analysis}, 249:113641, 2024.

\bibitem{surig2025gradient}
P.~S{\"u}rig.
\newblock Gradient estimates for {L}eibenson's equation on {R}iemannian
  manifolds.
\newblock {\em arXiv preprint arXiv:2506.07221}, 2025.

\bibitem{vazquez2015fundamental}
J.~L. V{\'a}zquez.
\newblock Fundamental solution and long time behavior of the porous medium
  equation in {H}yperbolic space.
\newblock {\em Journal de Math{\'e}matiques Pures et Appliqu{\'e}es},
  104(3):454--484, 2015.

\end{thebibliography}
	
	\emph{Universit\"{a}t Bielefeld, Fakult\"{a}t f\"{u}r Mathematik, Postfach
		100131, D-33501, Bielefeld, Germany}
	
	\texttt{philipp.suerig@uni-bielefeld.de}
\end{document}